\input amstex
\input amsppt.sty   
%
\catcode`\@=11

\def\East#1#2{\setboxz@h{$\m@th\ssize\;{#1}\;\;$}%
 \setbox@ne\hbox{$\m@th\ssize\;{#2}\;\;$}\setbox\tw@\hbox{$\m@th#2$}%
 \dimen@\minaw@
 \ifdim\wdz@>\dimen@ \dimen@\wdz@ \fi  \ifdim\wd@ne>\dimen@ \dimen@\wd@ne \fi
 \ifdim\wd\tw@>\z@
  \mathrel{\mathop{\hbox to\dimen@{\rightarrowfill}}\limits^{#1}_{#2}}%
 \else
  \mathrel{\mathop{\hbox to\dimen@{\rightarrowfill}}\limits^{#1}}%
 \fi}
\def\West#1#2{\setboxz@h{$\m@th\ssize\;\;{#1}\;$}%
 \setbox@ne\hbox{$\m@th\ssize\;\;{#2}\;$}\setbox\tw@\hbox{$\m@th#2$}%
 \dimen@\minaw@
 \ifdim\wdz@>\dimen@ \dimen@\wdz@ \fi \ifdim\wd@ne>\dimen@ \dimen@\wd@ne \fi
 \ifdim\wd\tw@>\z@
  \mathrel{\mathop{\hbox to\dimen@{\leftarrowfill}}\limits^{#1}_{#2}}%
 \else
  \mathrel{\mathop{\hbox to\dimen@{\leftarrowfill}}\limits^{#1}}%
 \fi}
\font\arrow@i=lams1
\font\arrow@ii=lams2
\font\arrow@iii=lams3
\font\arrow@iv=lams4
\font\arrow@v=lams5
\newbox\zer@
\newdimen\standardcgap
\standardcgap=40\p@
\newdimen\hunit
\hunit=\tw@\p@
\newdimen\standardrgap
\standardrgap=32\p@
\newdimen\vunit
\vunit=1.6\p@
\def\Cgaps#1{\RIfM@
  \standardcgap=#1\standardcgap\relax \hunit=#1\hunit\relax
 \else \nonmatherr@\Cgaps \fi}
\def\Rgaps#1{\RIfM@
  \standardrgap=#1\standardrgap\relax \vunit=#1\vunit\relax
 \else \nonmatherr@\Rgaps \fi}
\newdimen\getdim@
\def\getcgap@#1{\ifcase#1\or\getdim@\z@\else\getdim@\standardcgap\fi}
\def\getrgap@#1{\ifcase#1\getdim@\z@\else\getdim@\standardrgap\fi}
\def\cgaps#1{\RIfM@
 \cgaps@{#1}\edef\getcgap@##1{\i@=##1\relax\the\toks@}\toks@{}\else
 \nonmatherr@\cgaps\fi}
\def\rgaps#1{\RIfM@
 \rgaps@{#1}\edef\getrgap@##1{\i@=##1\relax\the\toks@}\toks@{}\else
 \nonmatherr@\rgaps\fi}
\def\Gaps@@{\gaps@@}
\def\cgaps@#1{\toks@{\ifcase\i@\or\getdim@=\z@}%
 \gaps@@\standardcgap#1;\gaps@@\gaps@@
 \edef\next@{\the\toks@\noexpand\else\noexpand\getdim@\noexpand\standardcgap
  \noexpand\fi}%
 \toks@=\expandafter{\next@}}
\def\rgaps@#1{\toks@{\ifcase\i@\getdim@=\z@}%
 \gaps@@\standardrgap#1;\gaps@@\gaps@@
 \edef\next@{\the\toks@\noexpand\else\noexpand\getdim@\noexpand\standardrgap
  \noexpand\fi}%
 \toks@=\expandafter{\next@}}
\def\gaps@@#1#2;#3{\mgaps@#1#2\mgaps@
 \edef\next@{\the\toks@\noexpand\or\noexpand\getdim@
  \noexpand#1\the\mgapstoks@@}%
 \global\toks@=\expandafter{\next@}%
 \DN@{#3}%
 \ifx\next@\Gaps@@\gdef\next@##1\gaps@@{}\else
  \gdef\next@{\gaps@@#1#3}\fi\next@}
\def\mgaps@#1{\let\mgapsnext@#1\FN@\mgaps@@}
\def\mgaps@@{\ifx\next\space@\DN@. {\FN@\mgaps@@}\else
 \DN@.{\FN@\mgaps@@@}\fi\next@.}
\def\mgaps@@@{\ifx\next\w\let\next@\mgaps@@@@\else
 \let\next@\mgaps@@@@@\fi\next@}
\newtoks\mgapstoks@@
\def\mgaps@@@@@#1\mgaps@{\getdim@\mgapsnext@\getdim@#1\getdim@
 \edef\next@{\noexpand\getdim@\the\getdim@}%
 \mgapstoks@@=\expandafter{\next@}}
\def\mgaps@@@@\w#1#2\mgaps@{\mgaps@@@@@#2\mgaps@
 \setbox\zer@\hbox{$\m@th\hskip15\p@\tsize@#1$}%
 \dimen@\wd\zer@
 \ifdim\dimen@>\getdim@ \getdim@\dimen@ \fi
 \edef\next@{\noexpand\getdim@\the\getdim@}%
 \mgapstoks@@=\expandafter{\next@}}
\def\changewidth#1#2{\setbox\zer@\hbox{$\m@th#2}%
 \hbox to\wd\zer@{\hss$\m@th#1$\hss}}
\atdef@({\FN@\ARROW@}
\def\ARROW@{\ifx\next)\let\next@\OPTIONS@\else
 \DN@{\csname\string @(\endcsname}\fi\next@}
\newif\ifoptions@
\def\OPTIONS@){\ifoptions@\let\next@\relax\else
 \DN@{\options@true\begingroup\optioncodes@}\fi\next@}
\newif\ifN@
\newif\ifE@
\newif\ifNESW@
\newif\ifH@
\newif\ifV@
\newif\ifHshort@
\expandafter\def\csname\string @(\endcsname #1,#2){%
 \ifoptions@\let\next@\endgroup\else\let\next@\relax\fi\next@
 \N@false\E@false\H@false\V@false\Hshort@false
 \ifnum#1>\z@\E@true\fi
 \ifnum#1=\z@\V@true\tX@false\tY@false\a@false\fi
 \ifnum#2>\z@\N@true\fi
 \ifnum#2=\z@\H@true\tX@false\tY@false\a@false\ifshort@\Hshort@true\fi\fi
 \NESW@false
 \ifN@\ifE@\NESW@true\fi\else\ifE@\else\NESW@true\fi\fi
 \arrow@{#1}{#2}%
 \global\options@false
 \global\scount@\z@\global\tcount@\z@\global\arrcount@\z@
 \global\s@false\global\sxdimen@\z@\global\sydimen@\z@
 \global\tX@false\global\tXdimen@i\z@\global\tXdimen@ii\z@
 \global\tY@false\global\tYdimen@i\z@\global\tYdimen@ii\z@
 \global\a@false\global\exacount@\z@
 \global\x@false\global\xdimen@\z@
 \global\X@false\global\Xdimen@\z@
 \global\y@false\global\ydimen@\z@
 \global\Y@false\global\Ydimen@\z@
 \global\p@false\global\pdimen@\z@
 \global\label@ifalse\global\label@iifalse
 \global\dl@ifalse\global\ldimen@i\z@
 \global\dl@iifalse\global\ldimen@ii\z@
 \global\short@false\global\unshort@false}
\newif\iflabel@i
\newif\iflabel@ii
\newcount\scount@
\newcount\tcount@
\newcount\arrcount@
\newif\ifs@
\newdimen\sxdimen@
\newdimen\sydimen@
\newif\iftX@
\newdimen\tXdimen@i
\newdimen\tXdimen@ii
\newif\iftY@
\newdimen\tYdimen@i
\newdimen\tYdimen@ii
\newif\ifa@
\newcount\exacount@
\newif\ifx@
\newdimen\xdimen@
\newif\ifX@
\newdimen\Xdimen@
\newif\ify@
\newdimen\ydimen@
\newif\ifY@
\newdimen\Ydimen@
\newif\ifp@
\newdimen\pdimen@
\newif\ifdl@i
\newif\ifdl@ii
\newdimen\ldimen@i
\newdimen\ldimen@ii
\newif\ifshort@
\newif\ifunshort@
\def\zero@#1{\ifnum\scount@=\z@
 \if#1e\global\scount@\m@ne\else
 \if#1t\global\scount@\tw@\else
 \if#1h\global\scount@\thr@@\else
 \if#1'\global\scount@6 \else
 \if#1`\global\scount@7 \else
 \if#1(\global\scount@8 \else
 \if#1)\global\scount@9 \else
 \if#1s\global\scount@12 \else
 \if#1H\global\scount@13 \else
 \Err@{\Invalid@@ option \string\0}\fi\fi\fi\fi\fi\fi\fi\fi\fi
 \fi}
\def\one@#1{\ifnum\tcount@=\z@
 \if#1e\global\tcount@\m@ne\else
 \if#1h\global\tcount@\tw@\else
 \if#1t\global\tcount@\thr@@\else
 \if#1'\global\tcount@4 \else
 \if#1`\global\tcount@5 \else
 \if#1(\global\tcount@10 \else
 \if#1)\global\tcount@11 \else
 \if#1s\global\tcount@12 \else
 \if#1H\global\tcount@13 \else
 \Err@{\Invalid@@ option \string\1}\fi\fi\fi\fi\fi\fi\fi\fi\fi
 \fi}
\def\a@#1{\ifnum\arrcount@=\z@
 \if#10\global\arrcount@\m@ne\else
 \if#1+\global\arrcount@\@ne\else
 \if#1-\global\arrcount@\tw@\else
 \if#1=\global\arrcount@\thr@@\else
 \Err@{\Invalid@@ option \string\a}\fi\fi\fi\fi
 \fi}
\def\ds@(#1;#2){\ifs@\else
 \global\s@true
 \sxdimen@\hunit \global\sxdimen@#1\sxdimen@\relax
 \sydimen@\vunit \global\sydimen@#2\sydimen@\relax
 \fi}
\def\dtX@(#1;#2){\iftX@\else
 \global\tX@true
 \tXdimen@i\hunit \global\tXdimen@i#1\tXdimen@i\relax
 \tXdimen@ii\vunit \global\tXdimen@ii#2\tXdimen@ii\relax
 \fi}
\def\dtY@(#1;#2){\iftY@\else
 \global\tY@true
 \tYdimen@i\hunit \global\tYdimen@i#1\tYdimen@i\relax
 \tYdimen@ii\vunit \global\tYdimen@ii#2\tYdimen@ii\relax
 \fi}
\def\da@#1{\ifa@\else\global\a@true\global\exacount@#1\relax\fi}
\def\dx@#1{\ifx@\else
 \global\x@true
 \xdimen@\hunit \global\xdimen@#1\xdimen@\relax
 \fi}
\def\dX@#1{\ifX@\else
 \global\X@true
 \Xdimen@\hunit \global\Xdimen@#1\Xdimen@\relax
 \fi}
\def\dy@#1{\ify@\else
 \global\y@true
 \ydimen@\vunit \global\ydimen@#1\ydimen@\relax
 \fi}
\def\dY@#1{\ifY@\else
 \global\Y@true
 \Ydimen@\vunit \global\Ydimen@#1\Ydimen@\relax
 \fi}
\def\p@@#1{\ifp@\else
 \global\p@true
 \pdimen@\hunit \divide\pdimen@\tw@ \global\pdimen@#1\pdimen@\relax
 \fi}
\def\L@#1{\iflabel@i\else
 \global\label@itrue \gdef\label@i{#1}%
 \fi}
\def\l@#1{\iflabel@ii\else
 \global\label@iitrue \gdef\label@ii{#1}%
 \fi}
\def\dL@#1{\ifdl@i\else
 \global\dl@itrue \ldimen@i\hunit \global\ldimen@i#1\ldimen@i\relax
 \fi}
\def\dl@#1{\ifdl@ii\else
 \global\dl@iitrue \ldimen@ii\hunit \global\ldimen@ii#1\ldimen@ii\relax
 \fi}
\def\s@{\ifunshort@\else\global\short@true\fi}
\def\uns@{\ifshort@\else\global\unshort@true\global\short@false\fi}
\def\optioncodes@{\let\0\zero@\let\1\one@\let\a\a@\let\ds\ds@\let\dtX\dtX@
 \let\dtY\dtY@\let\da\da@\let\dx\dx@\let\dX\dX@\let\dY\dY@\let\dy\dy@
 \let\p\p@@\let\L\L@\let\l\l@\let\dL\dL@\let\dl\dl@\let\s\s@\let\uns\uns@}
\def\slopes@{\\161\\152\\143\\134\\255\\126\\357\\238\\349\\45{10}\\56{11}%
 \\11{12}\\65{13}\\54{14}\\43{15}\\32{16}\\53{17}\\21{18}\\52{19}\\31{20}%
 \\41{21}\\51{22}\\61{23}}
\newcount\tan@i
\newcount\tan@ip
\newcount\tan@ii
\newcount\tan@iip
\newdimen\slope@i
\newdimen\slope@ip
\newdimen\slope@ii
\newdimen\slope@iip
\newcount\angcount@
\newcount\extracount@
\def\slope@{{\slope@i=\secondy@ \advance\slope@i-\firsty@
 \ifN@\else\multiply\slope@i\m@ne\fi
 \slope@ii=\secondx@ \advance\slope@ii-\firstx@
 \ifE@\else\multiply\slope@ii\m@ne\fi
 \ifdim\slope@ii<\z@
  \global\tan@i6 \global\tan@ii\@ne \global\angcount@23
 \else
  \dimen@\slope@i \multiply\dimen@6
  \ifdim\dimen@<\slope@ii
   \global\tan@i\@ne \global\tan@ii6 \global\angcount@\@ne
  \else
   \dimen@\slope@ii \multiply\dimen@6
   \ifdim\dimen@<\slope@i
    \global\tan@i6 \global\tan@ii\@ne \global\angcount@23
   \else
    \tan@ip\z@ \tan@iip \@ne
    \def\\##1##2##3{\global\angcount@=##3\relax
     \slope@ip\slope@i \slope@iip\slope@ii
     \multiply\slope@iip##1\relax \multiply\slope@ip##2\relax
     \ifdim\slope@iip<\slope@ip
      \global\tan@ip=##1\relax \global\tan@iip=##2\relax
     \else
      \global\tan@i=##1\relax \global\tan@ii=##2\relax
      \def\\####1####2####3{}%
     \fi}%
    \slopes@
    \slope@i=\secondy@ \advance\slope@i-\firsty@
    \ifN@\else\multiply\slope@i\m@ne\fi
    \multiply\slope@i\tan@ii \multiply\slope@i\tan@iip \multiply\slope@i\tw@
    \count@\tan@i \multiply\count@\tan@iip
    \extracount@\tan@ip \multiply\extracount@\tan@ii
    \advance\count@\extracount@
    \slope@ii=\secondx@ \advance\slope@ii-\firstx@
    \ifE@\else\multiply\slope@ii\m@ne\fi
    \multiply\slope@ii\count@
    \ifdim\slope@i<\slope@ii
     \global\tan@i=\tan@ip \global\tan@ii=\tan@iip
     \global\advance\angcount@\m@ne
    \fi
   \fi
  \fi
 \fi}%
}
\def\slope@a#1{{\def\\##1##2##3{\ifnum##3=#1\global\tan@i=##1\relax
 \global\tan@ii=##2\relax\fi}\slopes@}}
\newcount\i@
\newcount\j@
\newcount\colcount@
\newcount\Colcount@
\newcount\tcolcount@
\newdimen\rowht@
\newdimen\rowdp@
\newcount\rowcount@
\newcount\Rowcount@
\newcount\maxcolrow@
\newtoks\colwidthtoks@
\newtoks\Rowheighttoks@
\newtoks\Rowdepthtoks@
\newtoks\widthtoks@
\newtoks\Widthtoks@
\newtoks\heighttoks@
\newtoks\Heighttoks@
\newtoks\depthtoks@
\newtoks\Depthtoks@
\newif\iffirstnewCDcr@
\def\dotoks@i{%
 \global\widthtoks@=\expandafter{\the\widthtoks@\else\getdim@\z@\fi}%
 \global\heighttoks@=\expandafter{\the\heighttoks@\else\getdim@\z@\fi}%
 \global\depthtoks@=\expandafter{\the\depthtoks@\else\getdim@\z@\fi}}
\def\dotoks@ii{%
 \global\widthtoks@{\ifcase\j@}%
 \global\heighttoks@{\ifcase\j@}%
 \global\depthtoks@{\ifcase\j@}}
\def\prenewCD@#1\endnewCD{\setbox\zer@
 \vbox{%
  \def\arrow@##1##2{{}}%
  \rowcount@\m@ne \colcount@\z@ \Colcount@\z@
  \firstnewCDcr@true \toks@{}%
  \widthtoks@{\ifcase\j@}%
  \Widthtoks@{\ifcase\i@}%
  \heighttoks@{\ifcase\j@}%
  \Heighttoks@{\ifcase\i@}%
  \depthtoks@{\ifcase\j@}%
  \Depthtoks@{\ifcase\i@}%
  \Rowheighttoks@{\ifcase\i@}%
  \Rowdepthtoks@{\ifcase\i@}%
  \Let@
  \everycr{%
   \noalign{%
    \global\advance\rowcount@\@ne
    \ifnum\colcount@<\Colcount@
    \else
     \global\Colcount@=\colcount@ \global\maxcolrow@=\rowcount@
    \fi
    \global\colcount@\z@
    \iffirstnewCDcr@
     \global\firstnewCDcr@false
    \else
     \edef\next@{\the\Rowheighttoks@\noexpand\or\noexpand\getdim@\the\rowht@}%
      \global\Rowheighttoks@=\expandafter{\next@}%
     \edef\next@{\the\Rowdepthtoks@\noexpand\or\noexpand\getdim@\the\rowdp@}%
      \global\Rowdepthtoks@=\expandafter{\next@}%
     \global\rowht@\z@ \global\rowdp@\z@
     \dotoks@i
     \edef\next@{\the\Widthtoks@\noexpand\or\the\widthtoks@}%
      \global\Widthtoks@=\expandafter{\next@}%
     \edef\next@{\the\Heighttoks@\noexpand\or\the\heighttoks@}%
      \global\Heighttoks@=\expandafter{\next@}%
     \edef\next@{\the\Depthtoks@\noexpand\or\the\depthtoks@}%
      \global\Depthtoks@=\expandafter{\next@}%
     \dotoks@ii
    \fi}}%
  \tabskip\z@
  \halign{&\setbox\zer@\hbox{\vrule height10\p@ width\z@ depth\z@
   $\m@th\displaystyle{##}$}\copy\zer@
   \ifdim\ht\zer@>\rowht@ \global\rowht@\ht\zer@ \fi
   \ifdim\dp\zer@>\rowdp@ \global\rowdp@\dp\zer@ \fi
   \global\advance\colcount@\@ne
   \edef\next@{\the\widthtoks@\noexpand\or\noexpand\getdim@\the\wd\zer@}%
    \global\widthtoks@=\expandafter{\next@}%
   \edef\next@{\the\heighttoks@\noexpand\or\noexpand\getdim@\the\ht\zer@}%
    \global\heighttoks@=\expandafter{\next@}%
   \edef\next@{\the\depthtoks@\noexpand\or\noexpand\getdim@\the\dp\zer@}%
    \global\depthtoks@=\expandafter{\next@}%
   \cr#1\crcr}}%
 \Rowcount@=\rowcount@
 \global\Widthtoks@=\expandafter{\the\Widthtoks@\fi\relax}%
 \edef\Width@##1##2{\i@=##1\relax\j@=##2\relax\the\Widthtoks@}%
 \global\Heighttoks@=\expandafter{\the\Heighttoks@\fi\relax}%
 \edef\Height@##1##2{\i@=##1\relax\j@=##2\relax\the\Heighttoks@}%
 \global\Depthtoks@=\expandafter{\the\Depthtoks@\fi\relax}%
 \edef\Depth@##1##2{\i@=##1\relax\j@=##2\relax\the\Depthtoks@}%
 \edef\next@{\the\Rowheighttoks@\noexpand\fi\relax}%
 \global\Rowheighttoks@=\expandafter{\next@}%
 \edef\Rowheight@##1{\i@=##1\relax\the\Rowheighttoks@}%
 \edef\next@{\the\Rowdepthtoks@\noexpand\fi\relax}%
 \global\Rowdepthtoks@=\expandafter{\next@}%
 \edef\Rowdepth@##1{\i@=##1\relax\the\Rowdepthtoks@}%
 \colwidthtoks@{\fi}%
 \setbox\zer@\vbox{%
  \unvbox\zer@
  \count@\rowcount@
  \loop
   \unskip\unpenalty
   \setbox\zer@\lastbox
   \ifnum\count@>\maxcolrow@ \advance\count@\m@ne
   \repeat
  \hbox{%
   \unhbox\zer@
   \count@\z@
   \loop
    \unskip
    \setbox\zer@\lastbox
    \edef\next@{\noexpand\or\noexpand\getdim@\the\wd\zer@\the\colwidthtoks@}%
     \global\colwidthtoks@=\expandafter{\next@}%
    \advance\count@\@ne
    \ifnum\count@<\Colcount@
    \repeat}}%
 \edef\next@{\noexpand\ifcase\noexpand\i@\the\colwidthtoks@}%
  \global\colwidthtoks@=\expandafter{\next@}%
 \edef\Colwidth@##1{\i@=##1\relax\the\colwidthtoks@}%
 \colwidthtoks@{}\Rowheighttoks@{}\Rowdepthtoks@{}\widthtoks@{}%
 \Widthtoks@{}\heighttoks@{}\Heighttoks@{}\depthtoks@{}\Depthtoks@{}%
}
\newcount\xoff@
\newcount\yoff@
\newcount\endcount@
\newcount\rcount@
\newdimen\firstx@
\newdimen\firsty@
\newdimen\secondx@
\newdimen\secondy@
\newdimen\tocenter@
\newdimen\charht@
\newdimen\charwd@
\def\outside@{\Err@{This arrow points outside the \string\newCD}}
\newif\ifsvertex@
\newif\iftvertex@
\def\arrow@#1#2{\xoff@=#1\relax\yoff@=#2\relax
 \count@\rowcount@ \advance\count@-\yoff@
 \ifnum\count@<\@ne \outside@ \else \ifnum\count@>\Rowcount@ \outside@ \fi\fi
 \count@\colcount@ \advance\count@\xoff@
 \ifnum\count@<\@ne \outside@ \else \ifnum\count@>\Colcount@ \outside@\fi\fi
 \tcolcount@\colcount@ \advance\tcolcount@\xoff@
 \Width@\rowcount@\colcount@ \tocenter@=-\getdim@ \divide\tocenter@\tw@
 \ifdim\getdim@=\z@
  \firstx@\z@ \firsty@\mathaxis@ \svertex@true
 \else
  \svertex@false
  \ifHshort@
   \Colwidth@\colcount@
    \ifE@ \firstx@=.5\getdim@ \else \firstx@=-.5\getdim@ \fi
  \else
   \ifE@ \firstx@=\getdim@ \else \firstx@=-\getdim@ \fi
   \divide\firstx@\tw@
  \fi
  \ifE@
   \ifH@ \advance\firstx@\thr@@\p@ \else \advance\firstx@-\thr@@\p@ \fi
  \else
   \ifH@ \advance\firstx@-\thr@@\p@ \else \advance\firstx@\thr@@\p@ \fi
  \fi
  \ifN@
   \Height@\rowcount@\colcount@ \firsty@=\getdim@
   \ifV@ \advance\firsty@\thr@@\p@ \fi
  \else
   \ifV@
    \Depth@\rowcount@\colcount@ \firsty@=-\getdim@
    \advance\firsty@-\thr@@\p@
   \else
    \firsty@\z@
   \fi
  \fi
 \fi
 \ifV@
 \else
  \Colwidth@\colcount@
  \ifE@ \secondx@=\getdim@ \else \secondx@=-\getdim@ \fi
  \divide\secondx@\tw@
  \ifE@ \else \getcgap@\colcount@ \advance\secondx@-\getdim@ \fi
  \endcount@=\colcount@ \advance\endcount@\xoff@
  \count@=\colcount@
  \ifE@
   \advance\count@\@ne
   \loop
    \ifnum\count@<\endcount@
    \Colwidth@\count@ \advance\secondx@\getdim@
    \getcgap@\count@ \advance\secondx@\getdim@
    \advance\count@\@ne
    \repeat
  \else
   \advance\count@\m@ne
   \loop
    \ifnum\count@>\endcount@
    \Colwidth@\count@ \advance\secondx@-\getdim@
    \getcgap@\count@ \advance\secondx@-\getdim@
    \advance\count@\m@ne
    \repeat
  \fi
  \Colwidth@\count@ \divide\getdim@\tw@
  \ifHshort@
  \else
   \ifE@ \advance\secondx@\getdim@ \else \advance\secondx@-\getdim@ \fi
  \fi
  \ifE@ \getcgap@\count@ \advance\secondx@\getdim@ \fi
  \rcount@\rowcount@ \advance\rcount@-\yoff@
  \Width@\rcount@\count@ \divide\getdim@\tw@
  \tvertex@false
  \ifH@\ifdim\getdim@=\z@\tvertex@true\Hshort@false\fi\fi
  \ifHshort@
  \else
   \ifE@ \advance\secondx@-\getdim@ \else \advance\secondx@\getdim@ \fi
  \fi
  \iftvertex@
   \advance\secondx@.4\p@
  \else
   \ifE@ \advance\secondx@-\thr@@\p@ \else \advance\secondx@\thr@@\p@ \fi
  \fi
 \fi
 \ifH@
 \else
  \ifN@
   \Rowheight@\rowcount@ \secondy@\getdim@
  \else
   \Rowdepth@\rowcount@ \secondy@-\getdim@
   \getrgap@\rowcount@ \advance\secondy@-\getdim@
  \fi
  \endcount@=\rowcount@ \advance\endcount@-\yoff@
  \count@=\rowcount@
  \ifN@
   \advance\count@\m@ne
   \loop
    \ifnum\count@>\endcount@
    \Rowheight@\count@ \advance\secondy@\getdim@
    \Rowdepth@\count@ \advance\secondy@\getdim@
    \getrgap@\count@ \advance\secondy@\getdim@
    \advance\count@\m@ne
    \repeat
  \else
   \advance\count@\@ne
   \loop
    \ifnum\count@<\endcount@
    \Rowheight@\count@ \advance\secondy@-\getdim@
    \Rowdepth@\count@ \advance\secondy@-\getdim@
    \getrgap@\count@ \advance\secondy@-\getdim@
    \advance\count@\@ne
    \repeat
  \fi
  \tvertex@false
  \ifV@\Width@\count@\colcount@\ifdim\getdim@=\z@\tvertex@true\fi\fi
  \ifN@
   \getrgap@\count@ \advance\secondy@\getdim@
   \Rowdepth@\count@ \advance\secondy@\getdim@
   \iftvertex@
    \advance\secondy@\mathaxis@
   \else
    \Depth@\count@\tcolcount@ \advance\secondy@-\getdim@
    \advance\secondy@-\thr@@\p@
   \fi
  \else
   \Rowheight@\count@ \advance\secondy@-\getdim@
   \iftvertex@
    \advance\secondy@\mathaxis@
   \else
    \Height@\count@\tcolcount@ \advance\secondy@\getdim@
    \advance\secondy@\thr@@\p@
   \fi
  \fi
 \fi
 \ifV@\else\advance\firstx@\sxdimen@\fi
 \ifH@\else\advance\firsty@\sydimen@\fi
 \iftX@
  \advance\secondy@\tXdimen@ii
  \advance\secondx@\tXdimen@i
  \slope@
 \else
  \iftY@
   \advance\secondy@\tYdimen@ii
   \advance\secondx@\tYdimen@i
   \slope@
   \secondy@=\secondx@ \advance\secondy@-\firstx@
   \ifNESW@ \else \multiply\secondy@\m@ne \fi
   \multiply\secondy@\tan@i \divide\secondy@\tan@ii \advance\secondy@\firsty@
  \else
   \ifa@
    \slope@
    \ifNESW@ \global\advance\angcount@\exacount@ \else
      \global\advance\angcount@-\exacount@ \fi
    \ifnum\angcount@>23 \angcount@23 \fi
    \ifnum\angcount@<\@ne \angcount@\@ne \fi
    \slope@a\angcount@
    \ifY@
     \advance\secondy@\Ydimen@
    \else
     \ifX@
      \advance\secondx@\Xdimen@
      \dimen@\secondx@ \advance\dimen@-\firstx@
      \ifNESW@\else\multiply\dimen@\m@ne\fi
      \multiply\dimen@\tan@i \divide\dimen@\tan@ii
      \advance\dimen@\firsty@ \secondy@=\dimen@
     \fi
    \fi
   \else
    \ifH@\else\ifV@\else\slope@\fi\fi
   \fi
  \fi
 \fi
 \ifH@\else\ifV@\else\ifsvertex@\else
  \dimen@=6\p@ \multiply\dimen@\tan@ii
  \count@=\tan@i \advance\count@\tan@ii \divide\dimen@\count@
  \ifE@ \advance\firstx@\dimen@ \else \advance\firstx@-\dimen@ \fi
  \multiply\dimen@\tan@i \divide\dimen@\tan@ii
  \ifN@ \advance\firsty@\dimen@ \else \advance\firsty@-\dimen@ \fi
 \fi\fi\fi
 \ifp@
  \ifH@\else\ifV@\else
   \getcos@\pdimen@ \advance\firsty@\dimen@ \advance\secondy@\dimen@
   \ifNESW@ \advance\firstx@-\dimen@ii \else \advance\firstx@\dimen@ii \fi
  \fi\fi
 \fi
 \ifH@\else\ifV@\else
  \ifnum\tan@i>\tan@ii
   \charht@=10\p@ \charwd@=10\p@
   \multiply\charwd@\tan@ii \divide\charwd@\tan@i
  \else
   \charwd@=10\p@ \charht@=10\p@
   \divide\charht@\tan@ii \multiply\charht@\tan@i
  \fi
  \ifnum\tcount@=\thr@@
   \ifN@ \advance\secondy@-.3\charht@ \else\advance\secondy@.3\charht@ \fi
  \fi
  \ifnum\scount@=\tw@
   \ifE@ \advance\firstx@.3\charht@ \else \advance\firstx@-.3\charht@ \fi
  \fi
  \ifnum\tcount@=12
   \ifN@ \advance\secondy@-\charht@ \else \advance\secondy@\charht@ \fi
  \fi
  \iftY@
  \else
   \ifa@
    \ifX@
    \else
     \secondx@\secondy@ \advance\secondx@-\firsty@
     \ifNESW@\else\multiply\secondx@\m@ne\fi
     \multiply\secondx@\tan@ii \divide\secondx@\tan@i
     \advance\secondx@\firstx@
    \fi
   \fi
  \fi
 \fi\fi
 \ifH@\harrow@\else\ifV@\varrow@\else\arrow@@\fi\fi}
\newdimen\mathaxis@
\mathaxis@90\p@ \divide\mathaxis@36
\def\harrow@b{\ifE@\hskip\tocenter@\hskip\firstx@\fi}
\def\harrow@bb{\ifE@\hskip\xdimen@\else\hskip\Xdimen@\fi}
\def\harrow@e{\ifE@\else\hskip-\firstx@\hskip-\tocenter@\fi}
\def\harrow@ee{\ifE@\hskip-\Xdimen@\else\hskip-\xdimen@\fi}
\def\harrow@{\dimen@\secondx@\advance\dimen@-\firstx@
 \ifE@ \let\next@\rlap \else  \multiply\dimen@\m@ne \let\next@\llap \fi
 \next@{%
  \harrow@b
  \smash{\raise\pdimen@\hbox to\dimen@
   {\harrow@bb\arrow@ii
    \ifnum\arrcount@=\m@ne \else \ifnum\arrcount@=\thr@@ \else
     \ifE@
      \ifnum\scount@=\m@ne
      \else
       \ifcase\scount@\or\or\char118 \or\char117 \or\or\or\char119 \or
       \char120 \or\char121 \or\char122 \or\or\or\arrow@i\char125 \or
       \char117 \hskip\thr@@\p@\char117 \hskip-\thr@@\p@\fi
      \fi
     \else
      \ifnum\tcount@=\m@ne
      \else
       \ifcase\tcount@\char117 \or\or\char117 \or\char118 \or\char119 \or
       \char120\or\or\or\or\or\char121 \or\char122 \or\arrow@i\char125
       \or\char117 \hskip\thr@@\p@\char117 \hskip-\thr@@\p@\fi
      \fi
     \fi
    \fi\fi
    \dimen@\mathaxis@ \advance\dimen@.2\p@
    \dimen@ii\mathaxis@ \advance\dimen@ii-.2\p@
    \ifnum\arrcount@=\m@ne
     \let\leads@\null
    \else
     \ifcase\arrcount@
      \def\leads@{\hrule height\dimen@ depth-\dimen@ii}\or
      \def\leads@{\hrule height\dimen@ depth-\dimen@ii}\or
      \def\leads@{\hbox to10\p@{%
       \leaders\hrule height\dimen@ depth-\dimen@ii\hfil
       \hfil
      \leaders\hrule height\dimen@ depth-\dimen@ii\hskip\z@ plus2fil\relax
       \hfil
       \leaders\hrule height\dimen@ depth-\dimen@ii\hfil}}\or
     \def\leads@{\hbox{\hbox to10\p@{\dimen@\mathaxis@ \advance\dimen@1.2\p@
       \dimen@ii\dimen@ \advance\dimen@ii-.4\p@
       \leaders\hrule height\dimen@ depth-\dimen@ii\hfil}%
       \kern-10\p@
       \hbox to10\p@{\dimen@\mathaxis@ \advance\dimen@-1.2\p@
       \dimen@ii\dimen@ \advance\dimen@ii-.4\p@
       \leaders\hrule height\dimen@ depth-\dimen@ii\hfil}}}\fi
    \fi
    \cleaders\leads@\hfil
    \ifnum\arrcount@=\m@ne\else\ifnum\arrcount@=\thr@@\else
     \arrow@i
     \ifE@
      \ifnum\tcount@=\m@ne
      \else
       \ifcase\tcount@\char119 \or\or\char119 \or\char120 \or\char121 \or
       \char122 \or \or\or\or\or\char123\or\char124 \or
       \char125 \or\char119 \hskip-\thr@@\p@\char119 \hskip\thr@@\p@\fi
      \fi
     \else
      \ifcase\scount@\or\or\char120 \or\char119 \or\or\or\char121 \or\char122
      \or\char123 \or\char124 \or\or\or\char125 \or
      \char119 \hskip-\thr@@\p@\char119 \hskip\thr@@\p@\fi
     \fi
    \fi\fi
    \harrow@ee}}%
  \harrow@e}%
 \iflabel@i
  \dimen@ii\z@ \setbox\zer@\hbox{$\m@th\tsize@@\label@i$}%
  \ifnum\arrcount@=\m@ne
  \else
   \advance\dimen@ii\mathaxis@
   \advance\dimen@ii\dp\zer@ \advance\dimen@ii\tw@\p@
   \ifnum\arrcount@=\thr@@ \advance\dimen@ii\tw@\p@ \fi
  \fi
  \advance\dimen@ii\pdimen@
  \next@{\harrow@b\smash{\raise\dimen@ii\hbox to\dimen@
   {\harrow@bb\hskip\tw@\ldimen@i\hfil\box\zer@\hfil\harrow@ee}}\harrow@e}%
 \fi
 \iflabel@ii
  \ifnum\arrcount@=\m@ne
  \else
   \setbox\zer@\hbox{$\m@th\tsize@\label@ii$}%
   \dimen@ii-\ht\zer@ \advance\dimen@ii-\tw@\p@
   \ifnum\arrcount@=\thr@@ \advance\dimen@ii-\tw@\p@ \fi
   \advance\dimen@ii\mathaxis@ \advance\dimen@ii\pdimen@
   \next@{\harrow@b\smash{\raise\dimen@ii\hbox to\dimen@
    {\harrow@bb\hskip\tw@\ldimen@ii\hfil\box\zer@\hfil\harrow@ee}}\harrow@e}%
  \fi
 \fi}
\let\tsize@\tsize
\def\tsizenewCDlabels{\let\tsize@\tsize}
\def\ssizenewCDlabels{\let\tsize@\ssize}
\def\tsize@@{\ifnum\arrcount@=\m@ne\else\tsize@\fi}
\def\varrow@{\dimen@\secondy@ \advance\dimen@-\firsty@
 \ifN@ \else \multiply\dimen@\m@ne \fi
 \setbox\zer@\vbox to\dimen@
  {\ifN@ \vskip-\Ydimen@ \else \vskip\ydimen@ \fi
   \ifnum\arrcount@=\m@ne\else\ifnum\arrcount@=\thr@@\else
    \hbox{\arrow@iii
     \ifN@
      \ifnum\tcount@=\m@ne
      \else
       \ifcase\tcount@\char117 \or\or\char117 \or\char118 \or\char119 \or
       \char120 \or\or\or\or\or\char121 \or\char122 \or\char123 \or
       \vbox{\hbox{\char117 }\nointerlineskip\vskip\thr@@\p@
       \hbox{\char117 }\vskip-\thr@@\p@}\fi
      \fi
     \else
      \ifcase\scount@\or\or\char118 \or\char117 \or\or\or\char119 \or
      \char120 \or\char121 \or\char122 \or\or\or\char123 \or
      \vbox{\hbox{\char117 }\nointerlineskip\vskip\thr@@\p@
      \hbox{\char117 }\vskip-\thr@@\p@}\fi
     \fi}%
    \nointerlineskip
   \fi\fi
   \ifnum\arrcount@=\m@ne
    \let\leads@\null
   \else
    \ifcase\arrcount@\let\leads@\vrule\or\let\leads@\vrule\or
    \def\leads@{\vbox to10\p@{%
     \hrule height 1.67\p@ depth\z@ width.4\p@
     \vfil
     \hrule height 3.33\p@ depth\z@ width.4\p@
     \vfil
     \hrule height 1.67\p@ depth\z@ width.4\p@}}\or
    \def\leads@{\hbox{\vrule height\p@\hskip\tw@\p@\vrule}}\fi
   \fi
  \cleaders\leads@\vfill\nointerlineskip
   \ifnum\arrcount@=\m@ne\else\ifnum\arrcount@=\thr@@\else
    \hbox{\arrow@iv
     \ifN@
      \ifcase\scount@\or\or\char118 \or\char117 \or\or\or\char119 \or
      \char120 \or\char121 \or\char122 \or\or\or\arrow@iii\char123 \or
      \vbox{\hbox{\char117 }\nointerlineskip\vskip-\thr@@\p@
      \hbox{\char117 }\vskip\thr@@\p@}\fi
     \else
      \ifnum\tcount@=\m@ne
      \else
       \ifcase\tcount@\char117 \or\or\char117 \or\char118 \or\char119 \or
       \char120 \or\or\or\or\or\char121 \or\char122 \or\arrow@iii\char123 \or
       \vbox{\hbox{\char117 }\nointerlineskip\vskip-\thr@@\p@
       \hbox{\char117 }\vskip\thr@@\p@}\fi
      \fi
     \fi}%
   \fi\fi
   \ifN@\vskip\ydimen@\else\vskip-\Ydimen@\fi}%
 \ifN@
  \dimen@ii\firsty@
 \else
  \dimen@ii-\firsty@ \advance\dimen@ii\ht\zer@ \multiply\dimen@ii\m@ne
 \fi
 \rlap{\smash{\hskip\tocenter@ \hskip\pdimen@ \raise\dimen@ii \box\zer@}}%
 \iflabel@i
  \setbox\zer@\vbox to\dimen@{\vfil
   \hbox{$\m@th\tsize@@\label@i$}\vskip\tw@\ldimen@i\vfil}%
  \rlap{\smash{\hskip\tocenter@ \hskip\pdimen@
  \ifnum\arrcount@=\m@ne \let\next@\relax \else \let\next@\llap \fi
  \next@{\raise\dimen@ii\hbox{\ifnum\arrcount@=\m@ne \hskip-.5\wd\zer@ \fi
   \box\zer@ \ifnum\arrcount@=\m@ne \else \hskip\tw@\p@ \fi}}}}%
 \fi
 \iflabel@ii
  \ifnum\arrcount@=\m@ne
  \else
   \setbox\zer@\vbox to\dimen@{\vfil
    \hbox{$\m@th\tsize@\label@ii$}\vskip\tw@\ldimen@ii\vfil}%
   \rlap{\smash{\hskip\tocenter@ \hskip\pdimen@
   \rlap{\raise\dimen@ii\hbox{\ifnum\arrcount@=\thr@@ \hskip4.5\p@ \else
    \hskip2.5\p@ \fi\box\zer@}}}}%
  \fi
 \fi
}
\newdimen\goal@
\newdimen\shifted@
\newcount\Tcount@
\newcount\Scount@
\newbox\shaft@
\newcount\slcount@
\def\getcos@#1{%
 \ifnum\tan@i<\tan@ii
  \dimen@#1%
  \ifnum\slcount@<8 \count@9 \else \ifnum\slcount@<12 \count@8 \else
   \count@7 \fi\fi
  \multiply\dimen@\count@ \divide\dimen@10
  \dimen@ii\dimen@ \multiply\dimen@ii\tan@i \divide\dimen@ii\tan@ii
 \else
  \dimen@ii#1%
  \count@-\slcount@ \advance\count@24
  \ifnum\count@<8 \count@9 \else \ifnum\count@<12 \count@8
   \else\count@7 \fi\fi
  \multiply\dimen@ii\count@ \divide\dimen@ii10
  \dimen@\dimen@ii \multiply\dimen@\tan@ii \divide\dimen@\tan@i
 \fi}
\newdimen\adjust@
\def\Nnext@{\ifN@\let\next@\raise\else\let\next@\lower\fi}
\def\arrow@@{\slcount@\angcount@
 \ifNESW@
  \ifnum\angcount@<10
   \let\arrowfont@=\arrow@i \advance\angcount@\m@ne \multiply\angcount@13
  \else
   \ifnum\angcount@<19
    \let\arrowfont@=\arrow@ii \advance\angcount@-10 \multiply\angcount@13
   \else
    \let\arrowfont@=\arrow@iii \advance\angcount@-19 \multiply\angcount@13
  \fi\fi
  \Tcount@\angcount@
 \else
  \ifnum\angcount@<5
   \let\arrowfont@=\arrow@iii \advance\angcount@\m@ne \multiply\angcount@13
   \advance\angcount@65
  \else
   \ifnum\angcount@<14
    \let\arrowfont@=\arrow@iv \advance\angcount@-5 \multiply\angcount@13
   \else
    \ifnum\angcount@<23
     \let\arrowfont@=\arrow@v \advance\angcount@-14 \multiply\angcount@13
    \else
     \let\arrowfont@=\arrow@i \angcount@=117
  \fi\fi\fi
  \ifnum\angcount@=117 \Tcount@=115 \else\Tcount@\angcount@ \fi
 \fi
 \Scount@\Tcount@
 \ifE@
  \ifnum\tcount@=\z@ \advance\Tcount@\tw@ \else\ifnum\tcount@=13
   \advance\Tcount@\tw@ \else \advance\Tcount@\tcount@ \fi\fi
  \ifnum\scount@=\z@ \else \ifnum\scount@=13 \advance\Scount@\thr@@ \else
   \advance\Scount@\scount@ \fi\fi
 \else
  \ifcase\tcount@\advance\Tcount@\thr@@\or\or\advance\Tcount@\thr@@\or
  \advance\Tcount@\tw@\or\advance\Tcount@6 \or\advance\Tcount@7
  \or\or\or\or\or \advance\Tcount@8 \or\advance\Tcount@9 \or
  \advance\Tcount@12 \or\advance\Tcount@\thr@@\fi
  \ifcase\scount@\or\or\advance\Scount@\thr@@\or\advance\Scount@\tw@\or
  \or\or\advance\Scount@4 \or\advance\Scount@5 \or\advance\Scount@10
  \or\advance\Scount@11 \or\or\or\advance\Scount@12 \or\advance
  \Scount@\tw@\fi
 \fi
 \ifcase\arrcount@\or\or\advance\angcount@\@ne\else\fi
 \ifN@ \shifted@=\firsty@ \else\shifted@=-\firsty@ \fi
 \ifE@ \else\advance\shifted@\charht@ \fi
 \goal@=\secondy@ \advance\goal@-\firsty@
 \ifN@\else\multiply\goal@\m@ne\fi
 \setbox\shaft@\hbox{\arrowfont@\char\angcount@}%
 \ifnum\arrcount@=\thr@@
  \getcos@{1.5\p@}%
  \setbox\shaft@\hbox to\wd\shaft@{\arrowfont@
   \rlap{\hskip\dimen@ii
    \smash{\ifNESW@\let\next@\lower\else\let\next@\raise\fi
     \next@\dimen@\hbox{\arrowfont@\char\angcount@}}}%
   \rlap{\hskip-\dimen@ii
    \smash{\ifNESW@\let\next@\raise\else\let\next@\lower\fi
      \next@\dimen@\hbox{\arrowfont@\char\angcount@}}}\hfil}%
 \fi
 \rlap{\smash{\hskip\tocenter@\hskip\firstx@
  \ifnum\arrcount@=\m@ne
  \else
   \ifnum\arrcount@=\thr@@
   \else
    \ifnum\scount@=\m@ne
    \else
     \ifnum\scount@=\z@
     \else
      \setbox\zer@\hbox{\ifnum\angcount@=117 \arrow@v\else\arrowfont@\fi
       \char\Scount@}%
      \ifNESW@
       \ifnum\scount@=\tw@
        \dimen@=\shifted@ \advance\dimen@-\charht@
        \ifN@\hskip-\wd\zer@\fi
        \Nnext@
        \next@\dimen@\copy\zer@
        \ifN@\else\hskip-\wd\zer@\fi
       \else
        \Nnext@
        \ifN@\else\hskip-\wd\zer@\fi
        \next@\shifted@\copy\zer@
        \ifN@\hskip-\wd\zer@\fi
       \fi
       \ifnum\scount@=12
        \advance\shifted@\charht@ \advance\goal@-\charht@
        \ifN@ \hskip\wd\zer@ \else \hskip-\wd\zer@ \fi
       \fi
       \ifnum\scount@=13
        \getcos@{\thr@@\p@}%
        \ifN@ \hskip\dimen@ \else \hskip-\wd\zer@ \hskip-\dimen@ \fi
        \adjust@\shifted@ \advance\adjust@\dimen@ii
        \Nnext@
        \next@\adjust@\copy\zer@
        \ifN@ \hskip-\dimen@ \hskip-\wd\zer@ \else \hskip\dimen@ \fi
       \fi
      \else
       \ifN@\hskip-\wd\zer@\fi
       \ifnum\scount@=\tw@
        \ifN@ \hskip\wd\zer@ \else \hskip-\wd\zer@ \fi
        \dimen@=\shifted@ \advance\dimen@-\charht@
        \Nnext@
        \next@\dimen@\copy\zer@
        \ifN@\hskip-\wd\zer@\fi
       \else
        \Nnext@
        \next@\shifted@\copy\zer@
        \ifN@\else\hskip-\wd\zer@\fi
       \fi
       \ifnum\scount@=12
        \advance\shifted@\charht@ \advance\goal@-\charht@
        \ifN@ \hskip-\wd\zer@ \else \hskip\wd\zer@ \fi
       \fi
       \ifnum\scount@=13
        \getcos@{\thr@@\p@}%
        \ifN@ \hskip-\wd\zer@ \hskip-\dimen@ \else \hskip\dimen@ \fi
        \adjust@\shifted@ \advance\adjust@\dimen@ii
        \Nnext@
        \next@\adjust@\copy\zer@
        \ifN@ \hskip\dimen@ \else \hskip-\dimen@ \hskip-\wd\zer@ \fi
       \fi	
      \fi
  \fi\fi\fi\fi
  \ifnum\arrcount@=\m@ne
  \else
   \loop
    \ifdim\goal@>\charht@
    \ifE@\else\hskip-\charwd@\fi
    \Nnext@
    \next@\shifted@\copy\shaft@
    \ifE@\else\hskip-\charwd@\fi
    \advance\shifted@\charht@ \advance\goal@ -\charht@
    \repeat
   \ifdim\goal@>\z@
    \dimen@=\charht@ \advance\dimen@-\goal@
    \divide\dimen@\tan@i \multiply\dimen@\tan@ii
    \ifE@ \hskip-\dimen@ \else \hskip-\charwd@ \hskip\dimen@ \fi
    \adjust@=\shifted@ \advance\adjust@-\charht@ \advance\adjust@\goal@
    \Nnext@
    \next@\adjust@\copy\shaft@
    \ifE@ \else \hskip-\charwd@ \fi
   \else
    \adjust@=\shifted@ \advance\adjust@-\charht@
   \fi
  \fi
  \ifnum\arrcount@=\m@ne
  \else
   \ifnum\arrcount@=\thr@@
   \else
    \ifnum\tcount@=\m@ne
    \else
     \setbox\zer@
      \hbox{\ifnum\angcount@=117 \arrow@v\else\arrowfont@\fi\char\Tcount@}%
     \ifnum\tcount@=\thr@@
      \advance\adjust@\charht@
      \ifE@\else\ifN@\hskip-\charwd@\else\hskip-\wd\zer@\fi\fi
     \else
      \ifnum\tcount@=12
       \advance\adjust@\charht@
       \ifE@\else\ifN@\hskip-\charwd@\else\hskip-\wd\zer@\fi\fi
      \else
       \ifE@\hskip-\wd\zer@\fi
     \fi\fi
     \Nnext@
     \next@\adjust@\copy\zer@
     \ifnum\tcount@=13
      \hskip-\wd\zer@
      \getcos@{\thr@@\p@}%
      \ifE@\hskip-\dimen@ \else\hskip\dimen@ \fi
      \advance\adjust@-\dimen@ii
      \Nnext@
      \next@\adjust@\box\zer@
     \fi
  \fi\fi\fi}}%
 \iflabel@i
  \rlap{\hskip\tocenter@
  \dimen@\firstx@ \advance\dimen@\secondx@ \divide\dimen@\tw@
  \advance\dimen@\ldimen@i
  \dimen@ii\firsty@ \advance\dimen@ii\secondy@ \divide\dimen@ii\tw@
  \multiply\ldimen@i\tan@i \divide\ldimen@i\tan@ii
  \ifNESW@ \advance\dimen@ii\ldimen@i \else \advance\dimen@ii-\ldimen@i \fi
  \setbox\zer@\hbox{\ifNESW@\else\ifnum\arrcount@=\thr@@\hskip4\p@\else
   \hskip\tw@\p@\fi\fi
   $\m@th\tsize@@\label@i$\ifNESW@\ifnum\arrcount@=\thr@@\hskip4\p@\else
   \hskip\tw@\p@\fi\fi}%
  \ifnum\arrcount@=\m@ne
   \ifNESW@ \advance\dimen@.5\wd\zer@ \advance\dimen@\p@ \else
    \advance\dimen@-.5\wd\zer@ \advance\dimen@-\p@ \fi
   \advance\dimen@ii-.5\ht\zer@
  \else
   \advance\dimen@ii\dp\zer@
   \ifnum\slcount@<6 \advance\dimen@ii\tw@\p@ \fi
  \fi
  \hskip\dimen@
  \ifNESW@ \let\next@\llap \else\let\next@\rlap \fi
  \next@{\smash{\raise\dimen@ii\box\zer@}}}%
 \fi
 \iflabel@ii
  \ifnum\arrcount@=\m@ne
  \else
   \rlap{\hskip\tocenter@
   \dimen@\firstx@ \advance\dimen@\secondx@ \divide\dimen@\tw@
   \ifNESW@ \advance\dimen@\ldimen@ii \else \advance\dimen@-\ldimen@ii \fi
   \dimen@ii\firsty@ \advance\dimen@ii\secondy@ \divide\dimen@ii\tw@
   \multiply\ldimen@ii\tan@i \divide\ldimen@ii\tan@ii
   \advance\dimen@ii\ldimen@ii
   \setbox\zer@\hbox{\ifNESW@\ifnum\arrcount@=\thr@@\hskip4\p@\else
    \hskip\tw@\p@\fi\fi
    $\m@th\tsize@\label@ii$\ifNESW@\else\ifnum\arrcount@=\thr@@\hskip4\p@
    \else\hskip\tw@\p@\fi\fi}%
   \advance\dimen@ii-\ht\zer@
   \ifnum\slcount@<9 \advance\dimen@ii-\thr@@\p@ \fi
   \ifNESW@ \let\next@\rlap \else \let\next@\llap \fi
   \hskip\dimen@\next@{\smash{\raise\dimen@ii\box\zer@}}}%
  \fi
 \fi
}
\def\outnewCD@#1{\def#1{\Err@{\string#1 must not be used within \string\newCD}}}
\newskip\prenewCDskip@
\newskip\postnewCDskip@
\prenewCDskip@\z@
\postnewCDskip@\z@
\def\prenewCDspace#1{\RIfMIfI@
 \onlydmatherr@\prenewCDspace\else\advance\prenewCDskip@#1\relax\fi\else
 \onlydmatherr@\prenewCDspace\fi}
\def\postnewCDspace#1{\RIfMIfI@
 \onlydmatherr@\postnewCDspace\else\advance\postnewCDskip@#1\relax\fi\else
 \onlydmatherr@\postnewCDspace\fi}
\def\predisplayspace#1{\RIfMIfI@
 \onlydmatherr@\predisplayspace\else
 \advance\abovedisplayskip#1\relax
 \advance\abovedisplayshortskip#1\relax\fi
 \else\onlydmatherr@\prenewCDspace\fi}
\def\postdisplayspace#1{\RIfMIfI@
 \onlydmatherr@\postdisplayspace\else
 \advance\belowdisplayskip#1\relax
 \advance\belowdisplayshortskip#1\relax\fi
 \else\onlydmatherr@\postdisplayspace\fi}
\def\PrenewCDSpace#1{\global\prenewCDskip@#1\relax}
\def\PostnewCDSpace#1{\global\postnewCDskip@#1\relax}
\def\newCD#1\endnewCD{%
 \outnewCD@\cgaps\outnewCD@\rgaps\outnewCD@\Cgaps\outnewCD@\Rgaps
 \prenewCD@#1\endnewCD
 \advance\abovedisplayskip\prenewCDskip@
 \advance\abovedisplayshortskip\prenewCDskip@
 \advance\belowdisplayskip\postnewCDskip@
 \advance\belowdisplayshortskip\postnewCDskip@
 \vcenter{\vskip\prenewCDskip@ \Let@ \colcount@\@ne \rowcount@\z@
  \everycr{%
   \noalign{%
    \ifnum\rowcount@=\Rowcount@
    \else
     \global\nointerlineskip
     \getrgap@\rowcount@ \vskip\getdim@
     \global\advance\rowcount@\@ne \global\colcount@\@ne
    \fi}}%
  \tabskip\z@
  \halign{&\global\xoff@\z@ \global\yoff@\z@
   \getcgap@\colcount@ \hskip\getdim@
   \hfil\vrule height10\p@ width\z@ depth\z@
   $\m@th\displaystyle{##}$\hfil
   \global\advance\colcount@\@ne\cr
   #1\crcr}\vskip\postnewCDskip@}%
 \prenewCDskip@\z@\postnewCDskip@\z@
 \def\getcgap@##1{\ifcase##1\or\getdim@\z@\else\getdim@\standardcgap\fi}%
 \def\getrgap@##1{\ifcase##1\getdim@\z@\else\getdim@\standardrgap\fi}%
 \let\Width@\relax\let\Height@\relax\let\Depth@\relax\let\Rowheight@\relax
 \let\Rowdepth@\relax\let\Colwdith@\relax
}
\catcode`\@=\active
\mag=\magstep1
\hsize 30pc
\vsize 47pc
\def\nmb#1#2{#2}         
\def\cit#1#2{\ifx#1!\cite{#2}\else#2\fi} 
\def\totoc{}             
\def\idx{}               
\def\ign#1{}             
\redefine\o{\circ}

\define\al{\alpha}

\define\ga{\gamma}

\define\rh{\rho}
\define\si{\sigma}
\define\ta{\tau}

\define\ch{\chi}

\define\De{\Delta}

\define\Si{\Sigma}

\redefine\i{^{-1}}
\define\row#1#2#3{#1_{#2},\ldots,#1_{#3}}
\define\x{\times}
\define\pr{\operatorname{pr}}
\define\g{\frak g}
\define\Id{\operatorname{Id}}
\def\today{\ifcase\month\or
 January\or February\or March\or April\or May\or June\or
 July\or August\or September\or October\or November\or December\fi
 \space\number\day, \number\year}
\topmatter
\title Lifting smooth curves over invariants 
for representations of compact Lie groups
\endtitle
\author Dmitri Alekseevsky\\
Andreas Kriegl\\ Mark Losik\\ Peter W. Michor  
\endauthor
\date {November 24, 1997} \enddate
\leftheadtext{D\.V\. Alekseevsky, A\. Kriegl, M\. Losik, P\.W\. Michor}
\affil
Erwin Schr\"odinger International Institute of Mathematical Physics, 
Wien, Austria
\endaffil
\address 	 D\. V\. Alekseevsky: Center `Sophus Lie', 
Krasnokazarmennaya 6, 111250 Moscow, Russia
\endaddress
\address
A\. Kriegl, P\. W\. Michor: Institut f\"ur Mathematik, Universit\"at Wien,
Strudlhofgasse 4, A-1090 Wien, Austria
\endaddress
\email kriegl\@pap.univie.ac.at, Peter.Michor\@esi.ac.at \endemail
\address M\. Losik: Saratov State University, ul. Astrakhanskaya, 83, 410601 
Saratov, Russia \endaddress
\email losik\@scnit.saratov.su\endemail
\thanks  
P.W.M. was supported  
by `Fonds zur F\"orderung der wissenschaftlichen  
Forschung, Projekt P~10037~PHY'. 
\endthanks 
\abstract 
We show that one can lift locally real analytic curves from the orbit 
space of a compact Lie group representation, and that one can lift 
smooth curves even globally, but under an assumption.  
\endabstract
\endtopmatter

\document

\head Table of contents \endhead
\leftheadtext{\smc Alekseevsky, Kriegl, Losik, Michor}
\rightheadtext{\smc Lifting smooth curves over invariants}
\noindent 1. Introduction \leaders \hbox to 
1em{\hss .\hss }\hfill {\eightrm 1}\par  
\noindent 2. Invariants and their values \leaders \hbox to 
1em{\hss .\hss }\hfill {\eightrm 2}\par  
\noindent 3. Local lifting of curves over invariants 
\leaders \hbox to 1em{\hss .\hss }\hfill {\eightrm 5}\par  
\noindent 4. Global lifting of curves over invariants 
\leaders \hbox to 1em{\hss .\hss }\hfill {\eightrm 11}\par  

\head\totoc \nmb0{1}. Introduction
\endhead 

In \cit!{1} we investigated the following problem: 
Let
$$
P(t)=x^n-\si_1(t)x^{n-1}+\dots+(-1)^n\si_n(t) 
$$
be a polynomial with all roots real, smoothly parameterized by $t$ near 
0 in $\Bbb R$. 
Can we find $n$ smooth functions $x_1(t),\dots,x_n(t)$ of the 
parameter $t$ defined near 0, which are the roots of $P(t)$ for each $t$? 
We showed that 
this is possible under quite general conditions: real analyticity 
or no two roots should meet of infinite order. 

This problem can be reformulated in the following way: Let the 
symmetric group $S_n$ act $\Bbb R^n$ by permuting the coordinates 
(the roots), and consider the polynomial mapping 
$\si=(\si_1,\dots,\si_n):\Bbb R^n\to \Bbb R^n$ whose components are 
the elementary symmetric polynomials (the coefficients). Given a 
smooth curve $c:\Bbb R\to \si(\Bbb R^n)\subset \Bbb R^n$, is it possible 
to find a smooth lift $\bar c:\Bbb R\to \Bbb R^n$ with 
$\si\o\bar c=c$?

In this paper we tackle the following generalization of this problem. 
Consider an orthogonal representation of a compact Lie group $G$ 
on real vector space $V$ and
$\si=(\si_1,\dots,\si_n):V\to \Bbb R^n$, where $\si_1,\dots,\si_n$ is 
a system of generators for the algebra $\Bbb R[V]^G$ of invariant 
polynomials on $V$. 
Given an smooth curve $c:\Bbb R\to \Bbb R^m$ with values 
in $\si(V)$, does there exist a smooth lift $\bar c:\Bbb R\to V$ with 
$c=\si\o \bar c$? We show that this is true locally under similar 
conditions as for the case of the symmetric group and the elementary 
symmetric functions (i\.e\. the case of polynomials). Global results 
we get for smooth curves and for real analytic curves only in the 
case of polar representations (in particular for finite groups). 
In order to glue the 
local lifts of a real analytic curve we have to control their 
non-uniqueness. We try do to this by looking for orthogonal lifts 
(which are orthogonal to each orbit thy meet). Here many interesting 
unsolved problems arise, we suceed only in the case of polar 
representations. 

Similar lifting problems have been treated for smooth homotopies in 
\cit!{9}, and for diffeomorphisms (for Coxeter groups only) by 
\cit!{3}.
 
\head\totoc \nmb0{2}. Invariants and their values
\endhead 

\subhead\nmb.{2.1}. The setting \endsubhead
Let $G$ be a compact Lie group and $\rho:G\to O(V)$ an orthogonal 
representation in a real finite dimensional Euclidean vector space $V$ 
with an inner product $\langle \quad |\quad \rangle$.  
By a classical theorem of Hilbert and Nagata the 
algebra $\Bbb R[V]^G$ of invariant polynomials on $V$ is finitely 
generated. So let 
$\si_1,\dots,\si_n$ be a 
system of homogeneous generators of $\Bbb R[V]^G$ of degrees 
$d_1\leq\dots\leq d_n$.
Let $\si=(\si_1,\dots,\si_n):V\to \Bbb R^n$ 
and consider $\si(V)\subseteq \Bbb R^n$. 

\proclaim{Lemma} 
Let $(y_1,\dots,y_n)\in \si(V)\subseteq \Bbb R^n$. Then, for each 
$t\in \Bbb R$, $(t^{d_1}y_1,\dots,t^{d_n}y_n)\in \si(V)$. 
\endproclaim 

\demo{Proof} If $(y_1,\dots,y_n)=\si(v)$ for $v\in V$, then 
$(t^{d_1}y_1,\dots,t^{d_n}y_n)=\si(tv)$. 
\qed\enddemo 
 
Let $t\mapsto c(t)\in \si(V)$ be a smooth curve in $\Bbb R^n$, 
defined for $t$ near 0 in $\Bbb R$.
The curve $c$ is called \idx{\it (smoothly) liftable}
if there exists a smooth curve $t\mapsto \overline c(t)\in \Bbb R^n$ 
such that $c =\si\circ \overline{c}$. 
The curve  $\overline{c}$ is called a smooth lift of $c$. 

Similarly we consider real analytic lifts of real analytic curves.

\subhead\nmb.{2.2}. Orthogonal lifts \endsubhead 
A smooth lift $\overline c$ is called an \idx{\it orthogonal lift} if 
it is orthogonal to each orbit where it meets, i.e\. 
for each $t$ we have $\overline c'(t)\bot G.\overline c(t)$.
We suspect that an orthogonal lift of a curve which is nowhere 
infinitely flat exists and is unique up to the action of a constant in 
$G$. See \nmb!{4.2} below where we prove this for polar 
representations.

Note that orthogonal lifts may be not unique up to a constant in $G$ 
if we drop the assumption that $c$ is not infinitely flat. Namely, 
consider the action of $S^1$ by scalar multiplication on 
$\Bbb C\cong \Bbb R^2$. Then $\si(z)=|z|^2$. Any curve 
$c:\Bbb R\to \Bbb R$ which is non-negative and flat at a discrete set 
of $0$'s can be lifted to $\Bbb R^2$ along straight lines to 
$0\in \Bbb R^2$; at 0 we can change the angle arbitrily. 

Note that smooth lifts may exist but not orthogonal lifts, as example 
\cit!{1},~7.4 shows: there the compact group $O(2)$ acts by 
conjugation on the space $V$ of symmetric $2\x 2$-matrices, a polar 
representation. 
An orthogonal lift would lie in a section, without loss in the 
standard section of diagonal matrices, and this would mean a smooth choice 
of eigenvalues. 

\subhead\nmb.{2.3}. Removing fixed points \endsubhead
Consider the subspace 
$V^G:=\{v\in V: \rh(g)v=v\text{ for all }g\in G\}$ 
of $G$-invariant vectors of $V$ and the 
$G$-invariant decomposition $V=V^G\oplus (V^G)^{\bot}$. 
Obviously we have 
$\Bbb R[V]^G=\Bbb R[V^G]\otimes\Bbb R[(V^G)^{\bot}]^G$.  
 
Let us now order $\si_1,\dots,\si_n$ in a such way that 
$\si_1,\dots,\si_m$ $(m=\dim V^G)$ are coordinates on $V^G$ and 
$\si_{m+1},\dots,\si_n$ are generators of $\Bbb R[(V^G)^{\bot}]$. Put
$$
\si^{0}=(\si_1,\dots,\si_m):V^G\to \Bbb R^m 
\quad \text{and}\quad 
\si^{1}=(\si_{m+1},\dots,\si_n):(V^G)^{\bot}\to \Bbb R^{n-m}.
$$
Since $\si^{0}$ is an isomorphism of vector spaces and 
$\si=\si^{0}\times \si^{1}$, the study of $\si$ reduces to one of 
$\si^{1}$: a curve $c=(c^{0},c^{1})$ is (orthogonally) liftable over 
$\si=\si^0\x\si^1$ if and only if $c^1$ is (orthogonally) liftable 
over $\si^1$.

So without loss  
we suppose from now on that $V^G=\{0\}$. Then 
there exist no $G$-invariant linear forms on $V$ 
and $v\mapsto \langle v|v \rangle$ is 
a $G$-invariant homogeneous polynomial on $V$ of degree $2$; 
so we may also assume that $\si_1(v)=\langle v|v \rangle$ and $d_1=2$. 

If $y_1,\dots,y_n$ are the coordinates on $\Bbb R^n$, we have then:
\roster
\item For $x\in \si(V)$, $y_1(x)=0$ implies $x=0$. 
\endroster

\subhead\nmb.{2.4}. Description of $\si(V)$ \endsubhead
We recall some results of Procesi and Schwarz \cit!{4}. 
Let $I:=\{P\in \Bbb R[\Bbb R^n]: P\o \si=0\}$ 
be the ideal of relations of the generators $\si_1,\dots,\si_n$ in 
$\Bbb R[y_1,\dots,y_n]$, and let $Z\subset \Bbb R^n$ be the 
corresponding real algebraic set 
$Z:=\{y\in \Bbb R^n:P(y)=0\text{ for all }P\in I\}$.
Then $\si^*:\Bbb R[Z] := \Bbb R[\Bbb R^n]/I \to \Bbb R[V]^G$ 
is an isomorphism.
Since $\si$ is a polynomial map, its 
image $\si(V)\subset Z$ is a semialgebraic set. 
Since $G$ is compact, $\si$ is proper and separates orbits of $G$, 
thus it induces a homeomorphism between $V/G$ and $\si(V)$. 
 
Let $\langle \quad |\quad \rangle$ denote also the $G$-invariant 
dual inner product on $V^*$. 
The differentials $d\si_i:V\to V^*$ are 
$G$-equivariant, and the polynomials 
$v\mapsto \langle d\si_i(v)|d\si_j(v)\rangle$ are in $\Bbb R(V)^G$ 
and are entries of an $n\times n$ symmetric matrix valued polynomial
$$
B(v) := \pmatrix \langle d\si_1(v)|d\si_1(v)\rangle & \hdots & 
                    \langle d\si_1(v)|d\si_n(v)\rangle\\
                    \vdots &\ddots &\vdots\\
                 \langle d\si_n(v)|d\si_1(v)\rangle & \hdots &
                    \langle d\si_n(v)|d\si_n(v)\rangle
          \endpmatrix.
$$ 
There is a unique matrix valued polynomial
$\tilde B$ on $Z$ such that $B =\tilde B\o \si$.
 
For a real symmetric matrix $A$ let $A\ge 0$ indicate 
that $A$ is positive semidefinite. 
 
\proclaim{Theorem}(Procesi-Schwarz \cit!{4}) 
$\si(V)=\{z\in Z:\tilde B(z)\ge 0\}$. 
\endproclaim 
 
Denote by $x_i$ $(i=1,\dots,n)$ some orthonormal coordinates on $V$. 
Then $\si_1=\Sigma_{i=1}^nx_i^2$, $\si_i$ is a polynomial of degree $d_i$ 
in variables $x_i$ , and $d\si_i=\Sigma_{j=1}^n\frac{\partial 
\si_i}{\partial x_j}dx_j$. 
 
For each $1\le i_1 <\dots<i_s\le n$ $(s\le n)$, consider 
the matrix $B_{i_1,\dots,i_s}^{j_1,\dots,j_s}$ with entries 
$\langle d\si_{i_p}|d\si_{j_q}\rangle$ for  
$1\le p,q\le s$, a minor of $B$. 
We have the following equality: 
$$ 
B_{i_1,\dots,i_s}^{j_1,\dots,j_s}= 
\pmatrix \frac{\partial \si_{i_1}}{\partial x_1}&\hdots& 
          \frac{\partial \si_{i_1}}{\partial x_n}\\ 
     \vdots&\ddots&\vdots\\ 
     \frac{\partial \si_{i_s}}{\partial x_1}&\hdots&
          \frac{\partial \si_{i_s}}{\partial x_n}\endpmatrix 
\pmatrix \frac{\partial \si_{j_1}}{\partial x_1}&\hdots& 
          \frac{\partial \si_{j_1}}{\partial x_n}\\ 
     \vdots&\ddots&\vdots\\ 
     \frac{\partial \si_{j_s}}{\partial x_1}&\hdots&
          \frac{\partial \si_{j_s}}{\partial x_n}\endpmatrix^\top. 
$$ 
For the determinant 
$\De_{i_1,\dots,i_s}^{j_1,\dots,j_s}
     =\det B_{i_1,\dots,i_s}^{j_1,\dots,j_s}$ we have the following 
expression (see lemma \nmb!{2.5}): 
$$\align 
&\De_{i_1,\dots,i_s}^{j_1,\dots,j_s}= \sum_{1\le k_1<\dots,<k_s\le n}
D_{i_1,\dots,i_s}^{k_1,\dots,k_s}D_{j_1,\dots,j_s}^{k_1,\dots,k_s},\\
&\quad\text{ where }
D_{i_1,\dots,i_s}^{k_1,\dots,k_s} :=                     
\det\Bigl(\Bigl(\frac{\partial\si_{i_p}}{x_{k_q}}\Bigr)_{p,q=1,\dots,s}\Bigr)
\endalign$$ 
is the determinant of the minor of $d\si(v)$ 
determined by the rows with numbers $i_1,\dots,i_s$ and the columns with 
numbers $k_1,\dots,k_s$. 
Clearly, $\De_{i_1,\dots,i_s}^{j_1,\dots,j_s}(v)\neq 0$ means that 
the rank of $d\si(v)$ is larger or equal to $s$. 
The rank of $d\si(v)$ takes its maximum at regular $v$, and it 
is the codimension of any regular orbit in $V$.

A simple calculation shows that $\De_{i_1,\dots,i_s}^{j_1,\dots,j_s}$ is a 
polynomial  of degree
$d_{i_1}+ \dots+d_{i_s}+d_{j_1}+ \dots+d_{j_s}-2s$ on $V$. 
 
Since $\De_{i_1,\dots,i_s}^{j_1,\dots,j_s}$ is a $G$-invariant 
polynomial on $V$ there is a unique polynomial 
$\tilde\De_{i_1,\dots,i_s}^{j_1,\dots,j_s}$ on $Z$ such that 
$\De_{i_1,\dots,i_s}^{j_1,\dots,j_s}
     =\tilde\De_{i_1,\dots,i_s}^{j_1,\dots,j_s}\o \si$, by Schwarz' 
theorem \cit!{8}.

Then we have by Sylvester's theorem
$$\align
\si(V)&=\{z\in Z:\tilde B(z)\ge 0\}\\
&=\{z\in Z:\tilde \De_{i_1,\dots,i_k}^{i_1,\dots,i_k}(z)\ge 0 
     \text{ for }1\le i_1<\dots<i_k\le n\}. 
\endalign$$ 
Note also that 
$\De_1(v):=\De_1^1(v)=B_1^1(v)
     =\langle d\si_1(v),d\si_1(v) \rangle = 4\si_1(v)$ 
and consequently $\tilde \De_1(y)=4y_1$.

For further use let us also note here how removing of fixpoints as in 
\nmb!{2.3} affects the $\tilde\De_{i_1,\dots,i_s}^{j_1,\dots,j_s}$:
If $\si^0\x\si^1:V^G\x (V^G)^\bot\to \Bbb R^k\x \Bbb R^{n-k}$, then 
we have 
$\tilde\De(\si)_{1,2,\dots,k,i_1,\dots,i_s}^{1,2,\dots,k,j_1,\dots,j_s}
     =\tilde\De(\si^1)_{i_1,\dots,i_s}^{j_1,\dots,j_s}$.

\proclaim{\nmb.{2.5}. Lemma} For $s\le n$ let $A$ and $B$ be 
matrices with $s$ rows and $n$ columns. 

Then we have 
$$
\det(A.B^\top) = \sum_{1\le i_1<\dots<i_s\le n} 
\det A_{i_1,\dots,i_s}. \det B_{i_1,\dots,i_s},
$$
where $A_{i_1,\dots,i_s}$ is the quadratic minor of $A$ with columns 
${i_1,\dots,i_s}$. 
\endproclaim

\demo{Proof}
Let $a_i$ and $b_i$ for $i=1,\dots,s$ be the line 
with number $i$ of $A$ and $B$, 
respectively, so $a_i\in(\Bbb R^n)^*$ and $b_i^\top\in \Bbb R^n$. 
Then we have in terms of the standard basis $e_i$ of $\Bbb R^n$ and 
their dual basis $e^i$ of $(\Bbb R^n)^*$:
$$\align
&\det(AB^\top)=\det(\langle a_i,b_j^\top \rangle)
     =\langle a_1\wedge \dots\wedge a_s,
     b_1^\top\wedge \dots\wedge b_s^\top \rangle\\
&=\biggl\langle \sum_{i_1<\dots<i_s}\negmedspace
     \det(\langle a_p,e_{i_q}\rangle)
     e^{i_1}\wedge \dots\wedge e^{i_s}, 
     \negmedspace\sum_{j_1<\dots<j_s}\negmedspace
     \det(\langle e^{j_q},b_p^\top\rangle)
     e_{j_1}\wedge \dots\wedge e_{j_s}\biggr\rangle\\
&=\sum_{1\le i_1<\dots<i_s\le n}
     \det(\langle a_p,e_{i_q}\rangle)
     \det(\langle b_p,e_{i_q}\rangle)\\
&= \sum_{1\le i_1<\dots<i_s\le n} 
     \det A_{i_1,\dots,i_s}. \det B_{i_1,\dots,i_s}.
\quad\qed
\endalign$$
\enddemo

\proclaim{\nmb.{2.6}. Theorem} \cit!{7}
Let $\rh:G\to O(V)$ be an orthogonal, real, finite-dimensional 
representation of a compact Lie group $G$.
Let $\si_1,\dots,\si_n \in \Bbb R[V]^G$ be homogeneous generators for 
the algebra $\Bbb R[V]^G$ of invariant polynomials on $V$. For 
$v\in V$, let $N_v:=T_v(G.v)^\bot$ be the normal space to the 
orbit at $v$, and let $N_v^{G_v}$ be the subspace of those 
vectors which are invariant under the isotropy group $G_v$.

Then $\operatorname{grad}\si_1(v),\dots,\operatorname{grad}\si_n(v)$ 
span $N_v^{G_v}$ as a 
real vector space.
\endproclaim

\head\totoc\nmb0{3}. Local lifting of curves over invariants 
\endhead

\proclaim{\nmb.{3.1}. Lemma. Lifting at regular orbits}
Suppose that $c$ is a smooth curve in $\si(V)\subset \Bbb R^n$ such 
that $c(0)$, corresponds to a regular orbit.
Then $c$ is orthogonally smoothly liftable, locally near 0, and the 
orthogonal lift is uniquely determined by its initial value at $t=0$.
Thus any two orthogonal differentiable lifts differ by the action of a 
constant $g\in G$.

If $c$ is real analytic then any orthogonal lift is also real 
analytic.
\endproclaim

\demo{Proof}
Let $v\in V$ be a regular point with $\si(v)=c(0)$. 
Let $S$ be a normal slice at $v$, which can be chosen so small that 
the isotropy group $G_v$ acts trivial on $S$.
Then $G.S$ is open in $V$ and equivariantly diffeomorphic to 
$G/G_v\x S$ so that $S$ is mapped homeomorphically to an open 
neighborhood of the orbit through $v$ in the orbit space 
$V/G\cong \si(V)$. Thus $\si|S:S\to \si(V)\subset \Bbb R^n$ is 
injective. 

Then also $d(\si|S)_v:T_vS=T_v(G.v)^\bot\to \Bbb R^n$ is 
injective, by the following argument:
Choose coordinates $z_i$ on $S$ centered at $v$, let $h$ be a small 
bump function on $S$ which equals 1 near $v$. Since 
$G/G_v\x S\cong G.S$ equivariantly, $hz_i$ extends to a smooth 
$G$-invariant function on $V$, and by Schwarz' theorem \cit!{8} there 
exists a smooth function $f_i\in C^\infty(\Bbb R^n,\Bbb R)$ with 
$hz_i=f_i\o \si$. But then $z_i=d(hz_i)(v)=df_i(\si(v)).d\si(v)$, so 
that $d\si(v)$ is injective.

Now (making $S$ smaller if necessary) $\si|S: S\to \Bbb R^n$ is an 
embedding of a submanifold whose image is open in $\si(V)$, and by 
the implicit function theorem there are coordinates on $\Bbb R^n$ 
centered at $c(0)=\si(v)$ such that 
$\si|S:(z_1,\dots,z_m)\to (z_1,\dots,z_m,0,\dots,0)$. Then clearly 
$c$ is uniquely smoothly liftable near $t=0$ to $\tilde c:\Bbb R \to S$.
Via the embedding $S\subset V$ we get a smooth lift, which however is 
not orthogonal in general (only for polar representations).

Now $G.S\cong G/G_v\x S @>\pr_2>> S$ is a fiber bundle, and 
orthogonal projection onto the vertical bundle $\bigcup_{w}T(G.w)$ 
defines a (real analytic) (generalized or Ehresmann) connection on this 
fiber bundle, see \cit!{5}, section 9. The horizontal lift of the 
curve $\tilde c$ through any point in $G.v$ is then an orthogonal 
lift of $c$, and it is uniquely given by its initial value 
$\bar c(0)\in G.v$. 
\qed\enddemo

\subhead\nmb.{3.2}  \endsubhead 
For a smooth function $f$ defined near 0 in $\Bbb R$ let the 
\idx{\it multiplicity} or \idx{\it order of flatness} $m(f)$ at 0 be 
the supremum of all integer $p$ such that
$f(t)=t^pg(t)$ near $0$ for a smooth function $g$.
If $m(f)<\infty$ then $f(t)=t^{m(f)}g(t)$ where now $g(0)\ne 0$.
 
Let $c$ be a smooth curve in $V$ such that $m(c)\ge r>0$. 
Then $m(\si_i\o c)\ge rd_i$.

\proclaim{\nmb.{3.3}. Multiplicity Lemma} Let $c=(c_1,\dots,c_n)$ 
be a smooth curve in $\si(V)\subseteq \Bbb R^n$.
Then, for integers $r$, the following conditions are 
equivalent: 
\roster 
\item $m(c_i)\ge rd_i$ for $1\le i\le n$. 
\item $m(\tilde\De_{i_1,\dots,i_s}^{j_1,\dots,j_s}\o c)\ge 
     r(d_{i_1}+ \dots+d_{i_s}+d_{j_1}+ \dots+d_{j_s}-2s)$  
     for all $i_1\le\dots \le i_s$ and $i_j\le\dots \le j_s$. 
\item $m(c_1)\ge 2r$.  
\endroster 
\endproclaim 
\demo{Proof} 
(1) implies (2):  
By assumption $c_i(t)=t^{rd_i}b_i(t)$ for 
smooth $b_i$. Let us choose 
an arbitrary lift 
$\overline b$ of $b=(b_1,\dots,b_n)$, so that $\si\o \overline b=b$. 
Then 
$\si(t^r\overline b(t))=(t^{rd_1}b_1(t),\dots,t^{rd_n}b_n(t))=c(t)$, 
and we get
$$\align
\tilde \De_{i_1,\dots,i_s}^{j_1,\dots,j_s}(c(t))
&=(\tilde\De_{i_1,\dots,i_s}^{j_1,\dots,j_s}\o\si)(t^r\overline b(t)) 
=\De_{i_1,\dots,i_s}^{j_1,\dots,j_s}(t^r\overline b(t)) \\
&=t^{r(d_{i_1}+ \dots+d_{i_s}+d_{j_1}+ \dots+d_{j_s}-2s)}
\De_{i_1,\dots,i_s}^{j_1,\dots,j_s}(\overline b(t)) \\
&=t^{r(d_{i_1}+ \dots+d_{i_s}+d_{j_1}+ \dots+d_{j_s}-2s)}
\tilde\De_{i_1,\dots,i_s}^{j_1,\dots,j_s}(\si(\overline b(t))) \\
&=t^{r(d_{i_1}+ \dots+d_{i_s}+d_{j_1}+ \dots+d_{j_s}-2s)}
\tilde\De_{i_1,\dots,i_s}^{j_1,\dots,j_s}(b(t)) \\
\endalign$$
so 
$m(\tilde \De_{i_1,\dots,i_s}^{j_1,\dots,j_s}\o c)\ge 
     r(d_{i_1}+ \dots+d_{i_s}+d_{j_1}+ \dots+d_{j_s}-2s)$. 

(2) implies (3) since $\De_1=y_1$. 
 
(3) implies (1): Without loss let $r>0$. By 
\nmb!{2.3}.\therosteritem1 from $c_1(0)=0$ follows $c(0)=0$.
Therefore, $m(c_2),\dots, m(c_n)\geq 1$ and we have 
$c_k(t)=t^{m_k}c_{k,m_k}(t)$ for 
$k=1,\dots,n$, where the $m_k$ are positive integers 
and $c_{1,2r},c_{2,m_2},\dots,c_{n,m_n}$ 
are smooth functions, and where we may assume that either 
$m_k=m(c_k)<\infty$ or $m_k\ge rd_k$. 

Suppose now indirectly that for some $k> 1$ we have $m_k=m(c_k)<rd_k$.
Then we put 
$$ 
m:=\operatorname{min}(r,\frac{m_2}{d_2},\dots,\frac{m_n}{d_n}) < r.
$$ 
We consider the following continuous curve in $\Bbb R^n$ for 
$t\ge0$:
$$
c_{(m)}(t):=(t^{2r-2m}c_{1,2r}(t),t^{m_2-md_2}c_{2,m_2}(t),\dots, 
t^{m_n-md_n}c_{n,m_n}(t)). 
$$ 
By lemma \nmb!{2.1} the curve $c_{(m)}(t)\in \si(V)$ for $t>0$, and 
since $\si(V)$ is closed in $\Bbb R^n$ by \nmb!{2.4}, also 
$c_{(m)}(0)\in \si(V)$. Since $m<r$ the first coordinate of 
$c_{(m)}(t)$ vanishes at $t=0$, thus by 
\nmb!{2.3}.\therosteritem1 $c_{(m)}(0)=0$.
In particular for those $k$ with $m_k=md_k$ we have
$c_{k,m_k}(0)=0$ and, therefore, $m(c_k)>m_k$, a contradiction. 
\qed \enddemo 

\subhead\nmb.{3.4} \endsubhead
Let $r\in \Bbb N$ and let $c$ be a smooth curve in $\Bbb R^n$ with 
values in $\si(V)$ such that $c_k(t)=t^{m_k}c_{k,m_k}(t)$ for smooth 
$c_{k,m_k}$, and suppose that $m_k\ge rd_k$.
Consider the smooth curve 
$$ 
c_{(r)}(t)=(c_{1,m_1}(t)t^{m_1-rd_1},\dots,c_{n,m_n}(t)t^{m_n-rd_n}) 
\in \Bbb R^n.
$$ 
By Lemma \nmb!{2.1} $c_{(r)}(t)$ lies into $\si(V)$. 
If $c_{(r)}$ is liftable at $0$ and $\overline{c_{(r)}}$ is its 
smooth (real analytic) lift, 
$t\mapsto \bar c(t) :=t^r\overline{c_{(r)}}(t)$ is a smooth 
(real analytic) lift of $c$. 
If $\overline{c_{(r)}}$ is an orthogonal lift, then also 
$\bar c$, and conversely, since 
$\tfrac{d}{dt}(t^r\overline{c_{(r)}}(t))=rt^{r-1}\overline{c_{(r)}}(t)+ t^r\overline{c_{(r)}}'(t)$, 
where 
$\overline{c_{(r)}}(t)\bot |t^r.\overline{c_{(r)}}(t)|.S(V)
\supset t^r(G.\overline{c_{(r)}}(t))= G.(t^r\overline{c_{(r)}}(t))$. 
Moreover the orthogonal lift $\bar c$ is uniquely 
determined as orthogonal lift of $c$ up to the action of a constant 
element in $G$ if and only if $\overline{c_{(r)}}$ has this property as 
orthogonal lift of $c_{(r)}$.

Thus the problem of the existence of an orthogonal smooth lift of $c$ 
near $0$ and its uniqueness reduces to the one for 
$c_{(r)}$. 

\subhead\nmb.{3.5} \endsubhead
Let $v\ne0$ be a singular point of 
$V$. Since $V^G=0$ we have $G_v\neq G$. 
We consider the action of isotropy
group $G_v$ on the normal plane $N_v:= T_v(G.v)^\bot$, where we 
translate the vector space structure in such a way, that $v=0_{N_v}$. 
Since $v$ is a fixed point under $G_v$, this does not change the 
representation of $G_v$ on $N_v$. 
Let $\ta_1,\dots,\ta_m$ be a system of 
homogeneous generators of $\Bbb R[N_v]^{G_v}$ and 
$\ta=(\ta_1,\dots,\ta_m):N_v\to \Bbb R^m$ the corresponding map.

Since the restriction $\si_s|N_v\in\Bbb R[N_v]^{G_v}$ there exists a 
polynomial 
$l_s\in \Bbb R[\Bbb R^m]$ such that $\si_s|N_v=l_s(\ta_1,\dots,\ta_m)$.
Putting $l=(l_1,\dots,l_n):\Bbb R^m\to \Bbb R^n$ we get
$\si|N_v=l\circ \ta$.

\proclaim{Reduction Lemma}
In this situation let $c(t)=(c_1(t),\dots,c_n(t))$ be a smooth curve 
in $\si(V)$ such
that $c(0)=\si(v)$. Then there exists a smooth curve $c_v(t)$
in $\ta(N_v)$, defined near $t=0$, such that $c=l\circ c_v$. 
$$\rgaps{0.5;0.5}\newCD
N_v @()\L{\ta}@(1,0) @()\0(@(0,-2) & \Bbb R^m @()\L{l}@(0,-2) &  \\
 & & \Bbb R @()\L{c_v}@(-1,1) @()\l{c}@(-1,-1) \\
V @()\L{\si}@(1,0) & \Bbb R^n &  \\
\endnewCD$$
Moreover $c$ is locally smoothly liftable over $\si$ if and only if 
$c_v$ is locally smoothly liftable over $\ta$. 

Removing $G_v$-fixpoints as in \nmb!{2.3} then splits 
$$\ta=\ta^0\x \ta^1:N_v^{G_v}\x (N_v^{G_v})^\bot 
     \to \Bbb R^k\x \Bbb R^{m-k}$$
and we have $c_v=(c_v^0,c_v^1)$ where $c_v^0$ is smoothly liftable 
and where $c_v^1(0)=0$.

The real analytic analogon of this results is also true.
\endproclaim

\demo{Proof} Consider a normal slice $S$ at $v$, which is an open
neighborhood of $v$ in $N_v$, and the equivariant representation
$$\CD
G\x S @>{q}>> G\x_{G_v}S @>{\cong}>> G.S \\
@.               @VVV                    @VVrV \\
@.              G/G_v @>{\cong}>>        G.v, 
\endCD$$
for the $G$-invariant open neighborhood $G.S$ of the 
orbit $G.v$. Then the restriction homomorphism 
$\Bbb R[V]^G\to \Bbb R[N_v]^{G_v}$ is a monomorphism.

On the other hand, for each $f\in C^\infty(S)^{G_v}$ there exists
some $\tilde f\in C^\infty(V)^G$ such that $\tilde f|S$ agrees with 
$f$ near $v$, by the following argument: 
$(g,y)\mapsto f(y)$ is a smooth function $G\x S\to \Bbb R$ 
which is $G_v$-invariant and thus factors to a $G$-invariant function 
$G\x_{G_v}S\cong G.S\to \Bbb R$. 
We may multiply it by a 
$G$-invariant smooth function on $V$ with support in $G.S$ which 
equals $1$ near the orbit $G.v$ to obtain $\tilde f$. 

In particular we obtain $\tilde \ta_i\in C^\infty(V)^G$ with 
$\tilde\ta_i|S$ equals $\ta_i$ near $v$. 
By Schwarz' theorem there exist $h\in C^\infty(\Bbb R^n,\Bbb R^m)$
with $\tilde \ta=h\o\si$. 
Note that $\tilde \ta(G.S_1)=\ta(S_1)= h(\si(G.S_1))$.

Suppose $c=(c_1,\dots,c_n)$ is a smooth curve on $\si(V)$ such that
$c(0)=\si(v)$. Since the image of a slice under $\si$ is an open 
neighborhood of $\si(v)$ in $\si(V)$ we may 
suppose that $c(t)$ is contained in $\si(W)$, where
$W$ is a sufficiently small $G_v$-invariant neighborhood of $v$ in $N_v$.

On $W$ we have 
$l\o h\o \si = l\o \tilde \ta = l\o \ta = \si $, 
so $l\o h$ equals the identity on $\si(W)$, thus 
$c_v=h\o c$ 
is a smooth curve such that $c_v(0)=g(v)$ and 
$l\o c_v = c$. Since $h(\si(G.W))=\tilde\ta(G.W)=\ta(W)$ the curve 
$c_v$ has values in $\ta(N_v)$. 

Any lift $\overline{c_v}$ of $c_v$ over $\ta$ into $N_v\subset V$ is 
also a lift of $c$ over $\si$. 

Let conversely $\bar c$ be a smooth lift of $c$ 
over $\si$. Then $\bar c$ is locally in $G.S$; we may project it
to a smooth curve in the orbit $G.v\cong G/G_v$; this we can lift 
locally to a smooth curve $g$ in $G$, and then 
$g(t)\i.\bar c(t)\in S$ and describes a smooth lift of $c_v$ 
over $\ta$. This finishes the proof in the smooth case.

In the real analytic case we have to replace the smooth mapping 
$h\in C^\infty(\Bbb R^n,\Bbb R^m)$ from above by a real analytic one.
By omitting some of the $\ta_j$ if necessary let us first suppose 
without loss that the following two equivalent conditions hold:
\roster
\item"(a)" $\ta_1,\dots,\ta_m$ is a minimal system of homogeneous 
       generators of $\Bbb R[N_v]^{G_v}$, i\.e\. there does not exist 
       a polynomial relation of the form  
       $$\ta_j = P(\ta_1,\dots,\ta_{j-1},\ta_{j+1},\dots,\ta_m).$$
\item"(b)" $\ta_1,\dots,\ta_m$ induce an $\Bbb R$-basis of the vector 
       space $\Bbb R[N_v]_+^{G_v}/(\Bbb R[N_v]_+^{G_v})^2$, where we 
       denote by $\Bbb R[N_v]_+^{G_v}$ the space of all invariant 
       polynomials on  on $N_v$ which vanish at $0_{N_v}=v$.
\endroster
Clearly a polynomial relation as in (a) induces a nontrivial linear 
relation. If conversely (b) is wrong then we have 
$\sum a_i\ta_i = \sum_{|\al|\ge 2}c_\al ta^\al$, where the $\al$ are 
multiindices. Let $j$ be such that $a_j\ne0$ and that $\ta_j$ has 
minimal degree among them. Then we get 
$ta_j = \sum_{i>j}b_i\ta_i + \sum_{|\al|\ge 2}c'_\al ta^\al$.
If on the right hand side we omit all terms of degree different to 
$\deg \ta_j$ all terms involving $\ta_j$ vanish and we have a 
relation as forbidden in (a). So (a) and (b) are equivalent 
statements.

Note that we have 
$$\align
\si|N_v &= l \o \ta =: \si(v) + \tilde l \o \ta, \\
\ta &= h \o \si|N_v =: \tilde h\o (\si|N_v - \si(v)),\text{ near 
}v=0_{N_v},\\
&= \tilde h \o \tilde l\o \ta,
\endalign$$
where $\tilde l$ is polynomial and $\tilde h$ is smooth. 
Let us now pass to the infinite Taylor developments at 0, identify 
the polynomials with their Taylor series, and get
$$
\ta = j^\infty_0\tilde h \o \tilde l \o \ta = j^N_0 
\tilde h\o \tilde l \o \ta,
$$
since we may truncate the series of $h$ at some order higher than all 
the degrees of the $\ta_i$'s.
Since $\ta_1,\dots,\ta_m$ induces a real basis of 
$\Bbb R[N_v]_+^{G_v}/(\Bbb R[N_v]_+^{G_v})^2$ the last equation 
implies that for the derivative we have
$$
d(\tilde h\o \tilde l)(0) = d\tilde h(0)\o d\tilde l(0) = 
\Id_{\Bbb R^m}.
$$ 
Thus $\tilde l:\Bbb R^m\to \Bbb R^n$ is an immersion at $0_{N_v}=v$ 
and by the real analytic inverse function theorem 
there exists a real analytic mapping $H$ from an open 
neighborhood of 0 in $\Bbb R^n$ into $\Bbb R^m$ such that 
$H\o \tilde l = \Id_{\Bbb R^m}$ near 0. 
Using the real analytic mapping $H(\si(v)+\quad)$ instead of $h$ 
the foregoing proof for the smooth case can also be applied to the 
real analytic case.
\comment 

The following is inconclusive.
Suppose now that $\ga=\overline{c_v}$ is an orthogonal lift of $c_v$ 
in $N_v$, so that it is orthogonal to each $G_v$-orbit it meets. 
It is also a lift of $c$ over $\si$, but it need not be orthogonal as 
a curve into the $G$-space $V$.
We shall now look for a smooth curve $t\mapsto g(t)$ in $G$ with 
$g(0)=e$ such that $t\mapsto g(t).\ga(t)$ is orthogonal.

Let $\g$ be the Lie algebra of $G$, and let $\g_v=\{Y\in\g:Y.v=0\}$ 
be the Lie algebra of the isotropy group $G_v$. A smooth curve $g(t)$ 
with $g(0)=e$ is then uniquely determined by its \idx{\it left 
logarithmic derivative}, a smooth curve $t\mapsto X(t)$ in $\g$ with 
$g'(t)=g(t).X(t)$.
For such a $g(t)$ the following expression should vanish for all 
$Y\in \g$
$$\align
0&=\langle \partial_t (g(t).\ga(t)) , Y.g(t).\ga\rangle 
     =\langle g'(t).\ga(t)) + g(t).\ga'(t), Y.g(t).\ga\rangle\\
&=\langle g(t).X(t).\ga(t)) + g(t).\ga'(t), Y.g(t).\ga\rangle
     =\langle X(t).\ga(t)) + \ga'(t), g(t)\i.Y.g(t).\ga\rangle\\
\endalign$$
Let $\g=\g_v\oplus\frak m$, where $\frak m$ is a linear complement of 
$\g_v$ in $\g$. Let $Y_1,\dots,Y_k$ be a basis of $\frak m$. Then 
$Y_1.v,\dots,Y_k.v$ are linearly independent and orthogonal to $N_v$,
so for $w$ near $v$ in $N_v$ the vectors $Y_1.w,\dots,Y_k.w$ are 
still independent and and transversal to $N_v$. 

We now look for a smooth curve $X:\Bbb R\to \frak m$, defined near 
$t=0$, with $X(0)=0$, such that 
$t\mapsto\bar c(t):= \exp X(t).\overline{c_v}(t)$ is orthogonal in 
$V$. 
\endcomment
\qed\enddemo

\proclaim{\nmb.{3.6}. Theorem. Local real analytic lifts}
Let $\rh: G\to O(V)$ be an orthogonal representation of a compact Lie 
group $G$, and let $\si:V\to \Bbb R^n$ be the mapping whose 
components are generators of the algebra of invariant polynomials.
Let $c:\Bbb R\to \Bbb R^n$ be a real analytic curve with 
$c(\Bbb R)\subset \si(V)$.

Then there exists a real analytic lift
$\bar c$ in $V$ of $c$, locally near each $t\in \Bbb R$.
\endproclaim

\demo{Proof}
We may assume that $V^G=0$, and that all assumptions from \nmb!{2.3} 
hold.

We first show that there exist orthogonal lifts of $c$ locally near 
each point $t_0\in \Bbb R$, without loss $t_0=0$, 
through any $v\in\si\i(c(0))$. 
We do this by the following algorithm in 3 steps, 
which always stops, since it 
either gives an orthogonal local lift, 
or it reduces the lifting problem to a 
situation where the the group is replaced by an isotropy group of of 
nonzero vector, which either has smaller dimension, or at least less 
connected components (if it is already finite, e.g.). 

{\sl Step 1.}
If $c(0)\ne 0$ corresponds to a regular orbit, unique orthogonal real 
analytic lifts exist through all $v\in \si\i(c(0))$, by \nmb!{3.1}.

{\sl Step 2.}
If $c(0)\ne 0$ corresponds to a singular orbit, we consider the 
isotropy representation $G_v\to O(N_v)$ with invariant mapping 
$\ta:N_v\to \Bbb R^m$, we let $l:\Bbb R^m\to \Bbb R^n$ be a 
polynomial mapping with $\si|N_v = l\o \ta$. 
Then by the reduction lemma \nmb!{3.5} we may lift $c$ to 
$\tilde c:(\Bbb R,0)\to \Bbb R^m$ with $l\o \tilde c=c$, and it 
suffices to lift $\tilde c$ to $N_v$. Note that now $\tilde c(0)=0$, 
and that the group $G$ is replaced by the smaller (since $V^G=0$) 
group $G_v$, acting on $N_v$. By \nmb!{2.3} we may replace $N_v$ by 
$(N_v^{G_v})^\bot$, or assume that again 0 is the only fixpoint. 

{\sl Step 3.}
If $c(0)=0$ then $c_1(0)=0$ and $c_1=\tilde\De_1\o c\ge 0$, thus 
$m(c_1)= 2r$ for some integer $r\ge 1$ or $r=\infty$. In the latter 
case $c_1=0$ and by \nmb!{2.3}.\thetag1 $c=0$ is constant, which 
clearly can be lifted.

In the former case by the multiplicity lemma \nmb!{3.3} we have 
$m(c_i)\ge rd_i$ for some  
$r\ge 1$; so we can write $c_i(t)=t^{rd_i}c_{i,rd_i}(t)$ and we can 
consider the curve $c_{(r)}=(c_{1,2r},\dots c_{n,rd_n})$. By lemma 
\nmb!{2.1} this curve $c_{(r)}$ takes values in $\si(V)$, and if 
$\bar c_{(r)}$ is a lift of $c_{(r)}$, then 
$t\mapsto (t^r\bar c_{(r)}(t)$ is a real analytic lift of $c$. 
But now, since $m(c_1)=2r$, the first component $c_{1,2r}(0)$ of 
$c_{(r)}(0)$ does not vanish,
and we can continue in step 1 or step 2.
\qed\enddemo

\proclaim{\nmb.{3.7}. Theorem. Local smooth lifts}
Let $\rh: G\to O(V)$ be an orthogonal representation of a compact Lie 
group $G$, and let $\si:V\to \Bbb R^n$ be the mapping whose 
components are generators of the algebra of invariant polynomials.
Let $c:\Bbb R\to \Bbb R^n$ be a smooth curve with 
$c(\Bbb R)\subset \si(V)$ such that for some $r$ (not larger than the 
codimension of the regular orbits) we have:
\roster
\item $\tilde\De_{i_1,\dots,i_s}^{j_1,\dots,j_s}\o c=0$ whenever 
       $s>r$.
\item There exists a minor of size $r$ such that 
       $\tilde\De_{i_1,\dots,i_r}^{j_1,\dots,j_r}\o c$ does not 
       vanish of infinite order at $t_0$.
\endroster

Then there exists a smooth lift $\bar c$ in $V$ of $c$, locally near 
$t_0\in \Bbb R$.
\endproclaim

\demo{Proof}
We may assume that $V^G=0$, and that all assumptions from \nmb!{2.3} 
hold.

We first show that we can lift $c$ locally near  
$t_0\in \Bbb R$, without loss $t_0=0$, through any $v\in\si\i(c(0))$. 
We do this by the following algorithm in 3 steps, 
which always stops in step 1, since
step 2, followed by step 3, reduce the lifting problem to a 
situation where the the group is replaced by a smaller isotropy group of 
a nonzero vector and smaller $s$.

{\sl Step 1.}
If $c(0)\ne 0$ corresponds to a regular orbit, smooth lifts 
exist through all $v\in \si\i(c(0))$, by \nmb!{3.1}.

{\sl Step 2.}
If $c(0)\ne 0$ corresponds to a singular orbit, we consider the 
isotropy representation $G_v\to O(N_v)$ with invariant mapping 
$\ta:N_v\to \Bbb R^m$, we let $l:\Bbb R^m\to \Bbb R^n$ be a 
polynomial mapping with $\si|N_v = l\o \ta$. 
Then by the reduction lemma \nmb!{3.5} we may lift $c$ to 
$c_v:(\Bbb R,0)\to \Bbb R^m$ with $l\o c_v=c$, and it 
suffices to lift $c_v$ to $N_v$. Note that now $c_v(0)=0$, 
and that the group $G$ is replaced by the smaller (since $V^G=0$) 
group $G_v$, acting on $N_v$. 
Now we have to check that the assumption is still valid in this new 
situation. We have 
$$\align
\langle d\si_i(v),d\si_j(v) \rangle_V &= 
\langle d\si_i(v)|_{N_v},d\si_j(v)|_{N_v} \rangle_{N_v} \\
&= \biggl\langle \sum_k \frac{\partial l_i}{\partial z_k}(\ta(v)) d\ta_k(v),
   \sum_p \frac{\partial l_j}{\partial z_p}(\ta(v)) d\ta_p(v) 
   \biggr\rangle_{N_v} \\
&= \sum_k \frac{\partial l_i}{\partial z_k}(\ta(v)) 
     \sum_p \frac{\partial l_j}{\partial z_p}(\ta(v))
     \langle  d\ta_k(v),d\ta_p(v) \rangle_{N_v} \\
\endalign$$ 
so that by matrix theory (use lemma \nmb!{2.5} twice)
$$\align
(\tilde\De^G)_{i_1,\dots,i_s}^{j_1,\dots,j_s}\o c &= 
     \sum\Sb k_1<\dots<k_s\\m_1<\dots\le m_s\endSb 
     \det\biggl(\frac{\partial l_{i_p}}{\partial z_{k_q}}\o c_v\biggr)
     \biggl((\tilde\De^{G_v})_{k_1,\dots,k_s}^{m_1,\dots,m_s}\o c_v 
     \biggr)
     \det\biggl(\frac{\partial l_{j_p}}{\partial z_{m_q}}\o c_v\biggr).
\endalign$$
Now we remove all fixpoints (including $v$) in $N_v$, as in 
\nmb!{2.3}.  From the formula it follows that then for $s>r$ all 
$(\tilde\De^{G_v})_{1,2,\dots,m,i_1,\dots,i_{s-m}}
     ^{1,2,\dots,m,j_1,\dots,j_{s-m}}\o c_v$ 
vanish identically, 
whereas there is a minor of size $r$ such that   
$(\tilde\De^{G_v})_{1,2,\dots,m,i_1,\dots,i_{r-m}}
     ^{1,2,\dots,m,j_1,\dots,j_{r-m}}\o c_v$ 
does not vanish of infinite order at $0$. 
Note that the codimension of a regular 
orbit in $V$ equals the codimension of a regular $G_v$-orbit 
in the slice $S_v\subset N_v$. 
By the last remark in \nmb!{2.4} we have then 
$$
(\tilde\De(\ta^1))_{i_1,\dots,i_{s-k}}^{j_1,\dots,j_{s-k}}\o c_v^1 = 
  (\tilde\De^{G_v})^{1,2,\dots,k,j_1,\dots,j_{s-k}}
     _{1,2,\dots,k,i_1,\dots,i_{s-k}}\o c_v.$$
We now reproduced the starting situation for the curve 
$c_v^1$ in $\ta^1((N_v^{G_v})^\bot)$ with $c_v^1(0)=0$, and we 
reduced $r$ by $m=\dim (N_v^{G_v})\ge 1$. 

{\sl Step 3.}
If $c(0)=0$ then $c_1(0)=0$ and $c_1=\tilde\De_1\o c\ge 0$, thus 
$m(c_1)= 2r$ for some integer $r\ge 1$ or $r=\infty$. In the latter 
case the multiplicity lemma \nmb!{3.3} gives 
$m(\tilde\De_{i_1,\dots,i_s}^{j_1,\dots,j_s}\o c)=\infty$, a 
contradiction. 

In the former case by the multiplicity lemma \nmb!{3.3} we have 
$m(c_i)\ge rd_i$ for some  
$r\ge 1$; so we can write $c_i(t)=t^{rd_i}c_{i,rd_i}(t)$ and we can 
consider the curve $c_{(r)}=(c_{1,2r},\dots c_{n,rd_n})$. By lemma 
\nmb!{2.1} this curve $c_{(r)}$ takes values in $\si(V)$, and if 
$\bar c_{(r)}$ is a lift of $c_{(r)}$, then 
$t\mapsto (t^r\bar c_{(r)}(t)$ is a smooth lift of $c$. 
But now, since $m(c_1)=2r$, the first component $c_{1,2r}(0)$ of 
$c_{(r)}(0)$ does not vanish, and (see the proof of \nmb!{3.3})
$$\tilde\De_{i_1,\dots,i_s}^{j_1,\dots,j_s}\o c_{(r)}=
t^{-r(d_{i_1}+ \dots+d_{i_s}+d_{j_1}+ \dots+d_{j_s}-2s)}
(\tilde\De_{i_1,\dots,i_s}^{j_1,\dots,j_s}\o c)$$ 
is not infinitely flat.
So we can continue in step 1 or step 2.

This finishes the proof of the existence of local lifts.
\qed\enddemo

\proclaim{\nmb.{3.8}. Lemma}
Let $\rh: G\to O(V)$ be an orthogonal representation of a compact Lie 
group $G$, and let $\si:V\to \Bbb R^n$ be the mapping whose 
components are generators of the algebra of invariant polynomials.
Let $c:\Bbb R\to \Bbb R^n$ be a smooth curve with 
$c(\Bbb R)\subset \si(V)$ which is nowhere infinitely flat.
Suppose that $\bar c_1, \bar c_2: I\to V$ are smooth lifts of $c$ on 
an open interval $I$. 

Then for each $t_0\in I$ there exists a smooth curve 
$g$ in $G$ defined near $t_0$  
such that $\bar c_1(t)=g(t).\bar c_2(t)$ for all $t$ near $t_0$. 
The real analytic version of this result is also true.
\endproclaim

\demo{Proof}
We prove this by induction on the size 
(dimension, and number of connected components in the case of the 
same dimension)
of $G$. 

Without loss let $t_0=0$. We also choose $g_0\in G$ with 
$\bar c_1(0)=g_0.\bar c_2(0)$, or assume without loss that 
$\bar c_1(0)=\bar c_2(0)$. 

Let us first remove the fixpoints $V^G$ and let us assume that 
the assumptions from \nmb!{2.3} hold. Note that then $c(0)=0$ if and 
only if $\bar c_i(0)=0$.

{\sl Step 1.} 
If $c(0)=0$ then $\bar c_i(t)=t^r\bar c_{i,r}(t)$ with 
$\bar c_{i,r}(0)\ne 0$, and $\si\o \bar c_{i,r} = c_{(r)}$ in the 
setting of step 3 of \nmb!{3.7}. If we can find $g(t)\in G$ with 
$g(t).\bar c_{1,r}(t)=\bar c_{2,r}(t)$ then we also have
$g(t).\bar c_1(t)= t^r.g(t).\bar c_{1,r}(t)= t^r.\bar c_{2,r}(t)= 
\bar c_2(t)$. Thus we may assume that $c(0)\ne 0$.

{\sl Step 2.}
If $c(0)\ne 0$, then for a normal slice $S_v$ at 
$v=\bar c_1(0)$ we know that 
$r:G.S_v\cong G\x_{G_v} S_v\to G/G_v\cong G.v$ is the 
projection of a fiber bundle associated to the principal bundle 
$G\to G/G_v$. Then $r\o \bar c_1$ and $r\o \bar c_2$ are two smooth 
curves in $G/G_v$ defined near $t=0$, which admit smooth lifts $g_1$ 
and $g_2$ into $G$ (via the horizontal lift of a principal 
connection, say). Then $t\mapsto g_j(t)\i.\bar c_j(t)$ are two smooth 
curves in $S_v$, lifts of $c$. We may now remove the fixpoints 
$N_v^{G_v}$ containing $v$, 
then we have to apply step 1, and have reduced the 
situation to the smaller group $G_v$.
Note that if $v$ is a regular point then $G_v$ acts trivially on 
$N_v$ and the two lifts are automatically the same.

In the real analytic situation the proof is the same: one has to use 
a real analytic principal connection in step 2.
\qed\enddemo

\head\totoc\nmb0{4}. Global lifting of curves over invariants \endhead

\proclaim{\nmb.{4.1}. Theorem. Global smooth lifts}
Let $\rh: G\to O(V)$ be an orthogonal representation of a compact Lie 
group $G$, and let $\si:V\to \Bbb R^n$ be the mapping whose 
components are generators of the algebra of invariant polynomials.
Let $c:\Bbb R\to \Bbb R^n$ be a smooth curve with 
$c(\Bbb R)\subset \si(V)$ 
such that for some $r$ (not larger than the 
codimension of the regular orbits) we have:
\roster
\item $\tilde\De_{i_1,\dots,i_s}^{j_1,\dots,j_s}\o c=0$ for any minor 
       of size $s>r$.
\item For each $t$ there exists a minor of size $r$ such that 
       $\tilde\De_{i_1,\dots,i_r}^{j_1,\dots,j_r}\o c$ does not 
       vanish of infinite order at $t$. 
\endroster
Then there exists a global smooth lift $\bar c: \Bbb R\to V$ with 
$\si\o \bar c= c$. 
\endproclaim

\demo{Proof}
By \nmb!{3.7} there exist local smooth lifts near any $t\in \Bbb R$.
Now let $I$ be a maximal open interval such that a smooth lift 
$\bar c_1: I \to V$ of $c$ exists. We claim that $I=\Bbb R$. If not, 
it is bounded from above (say), let $t_0$ be its upper boundary 
point. By the first part of the proof there exists a local smooth 
lift $\bar c_2$ of $c$ near $t_0$, and a $t_1<t_0$ such that both 
$\bar c_1$ and $\bar c_2$ are defined near $t_1$. By lemma 
\nmb!{3.8} there exists a smooth curve $g$ in $G$, locally defined 
near $t_1$, such that $\bar c_1(t)=g(t).\bar c_2(t)$. We consider the 
right logarithmic derivative $X(t)=g'(t).g(t)\i\in \g$ and choose a 
smooth function $\ch(t)$ which is 1 for $t\le t_1$ and becomes 0 
before $g$ ceases to exist. Then $Y(t)=\ch(t)X(t)$ is smooth and 
defined near $[t_1,\infty)$. The differential equation 
$h'(t)=Y(t).h(t)$ with initial condition $h(t_1)=g(t_1)$ then has a 
solution $h$ in $G$ defined near $[t_1,\infty)$ which coincides with 
$g$ below $t_1$. Then
$$
\bar c(t):=\cases \bar c_1(t) &\text{ for }t\le t_1\\
                  h(t).\bar c_2(t) &\text{ for }t\ge t_1 \endcases
$$
is a smooth lift of $c$ on a larger interval.
\qed\enddemo

\proclaim{\nmb.{4.2}. Theorem. Polar representations}
Let $\rh: G\to O(V)$ be a polar orthogonal representation of a compact Lie 
group $G$, see \cit!{2}, 
and let $\si:V\to \Bbb R^n$ be the mapping whose 
components are generators of the algebra of invariant polynomials.
Let $c:\Bbb R\to \Bbb R^n$ be a curve with 
$c(\Bbb R)\subset \si(V)$ which is either real analytic, or smooth 
but for some $r$ satisfies the two conditions 
\nmb!{4.1}.\therosteritem1 and \therosteritem2. 

Then there exists a global orthogonal real analytic or smooth lift $\bar c: 
\Bbb R\to V$ with $\si\o \bar c= c$.
If moreover $G$ is connected then an orthogonal lift is 
unique up to the action of a constant in $G$.
\endproclaim

A representation is polar if there exists a linear subspace 
$\Si\subset V$, called a {\it section} or a \idx{\it Cartan 
subspace}, which meets each orbit orthogonally. See \cit!{2} and 
\cit!{6}. Then through each point of $V$ there is a section, which 
is the normal to the tangent of any regular orbit it meets. The trace 
of the $G$-action is the action of the \idx{\it generalized Weyl group} 
$W(\Si)=N_G(\Si)/Z_G(\Si)$ on $\Si$, which is a finite group, 
and is a reflection group for connected $G$. We shall also need the 
following generalization of the Chevalley restriction theorem, which 
is due to Chuu Lian Terng, see \cit!{6}, 4.12, or \cit!{10}, 
theorem~D: 
{\it For a polar representation the algebra $\Bbb R[V]^G$ of 
$G$-invariant polynomials on $V$ is isomorphic to the algebra 
$\Bbb R[\Si]^{W(\Si)}$ of Weyl group-invarint polynomials on $\Si$, 
via restriction.}

\demo{Proof} Let $\Si$ be a section.
By Terng's result, $\si|\Si:\Si\to \Bbb R^n$ is an 
invariant mapping for the representation $W=W(\Si)\to O(\Si)$ whose 
components are generators for the algebra of all invariant 
polynomials. If $c$ is a smooth curve with properties 
\nmb!{4.1}.\therosteritem1 and \therosteritem2, then by theorem 
\nmb!{4.1} there exists a global lift $\bar c:\Bbb R\to \Si$, which
as curve in $V$ is orthogonal to each $G$-orbit it meets, by the 
properties of $\Si$. Note for further use that $\bar c$ is nowhere 
infinitely flat, since otherwise property \nmb!{4.1}.\therosteritem2 
would not hold for $c$.

If $c$ is real analytic there are local lifts 
over $\si|\Si$
into $\Si$ by theorem \nmb!{3.6}. 

We claim that these local lifts are 
unique up to the action of a constant element in $W$. Namely, let 
$\bar c_1$ and $\bar c_2$ be real analytic lifts defined on an 
interval $I$. Choose a convergent sequence $t_i\in I$ and elements 
$\al_i\in W$ with $\al_i.\bar c_1(t_i)=\bar c_2(t_i)$. Since $W$ is 
finite, by passing to a subsequence we may assume that all 
$\al_i=\al\in W$. But then the real analytic curves $\bar c_2$ and 
$\al.\bar c_1$ coincide on a converging sequence, so they coincide 
on the whole interval. 

Thus we may glue the local lifts to a global real analytic lift 
$\bar c$ in $\Si$, which as curve in $V$ is then orthogonal. 

It remains to show that for two orthogonal lifts 
$\bar c_1,\bar c_2:\Bbb R\to V$ of $c$ 
there is a constant element $g\in G$ with 
$g.\bar c_1(t)=\bar c_2(t)$ for all $t$, 
or that they are equal if 
$\bar c_1(0)=\bar c_2(0)$. We may also assume that $\bar c_1$ lies in 
a section $\Si$, by the first assertion. 

 From lemma \nmb!{3.8} we get that  $\bar c_1(t)=g(t).\bar c_2(t)$ for 
some smooth or real analytic curve $g:I\to G$, locally near each 
$t_0$, with $g(t_0)=e$. We 
consider the right  
logarithmic derivative $X(t):= g'(t).g(t)\i\in \g$. Differentiating 
$\bar c_1(t)=g(t).\bar c_2(t)$ we get
$$
\bar c_1'(t)-g(t).\bar c_2'(t) = X(t).g(t).\bar c_2(t) = 
X(t).\bar c_1(t),
$$
where the left hand side is orthogonal to the orbit through 
$\bar c_1(t)$, and the right hand side is tangential to it,
so both sides are zero and $X(t)$ lies in the isotropy Lie algebra 
$\g_{\bar c_1(t)}$ for each $t$.  

The curve $\bar c_1$ lies in the section $\Si$, and $\Si$ carries the 
structure of walls and Weyl chambers from the action of the Coxeter 
group $W$, which encodes the stratification into different orbit 
types of the $G$-action on $V$. In particular the isotropy algebras 
$\g_v$ are constant for $v$ in each submanifold of the stratification of 
$\Si$ given by intersections of walls. $\g_v$ becomes larger if 
$v$ is contained in more walls. 

For the final argument we first assume that $\bar c_1$ is real analytic. 
Then it hits intersections of walls 
only in isolated points if it is not contained in them forever. But 
this means, that $X(t)\in\g_{\bar c_1(t)}$ lies in a fixed isotropy 
algebra (the smallest which it meets) 
for a dense open subset of $t's$, thus for all $t$.
But then $g(t)$ lies in the corresponding isotropy group and thus 
does not move $\bar c_1$. So $\bar c_2=\bar c_1$. 

If $\bar c_1$ is not real analytic we have to use a more complicated 
argument to get a dense set of $t$'s as above.
We consider the stratification of $\Si$ and let $L$ be the 
intersection of all walls containing $\bar c_1(\Bbb R)$. Let 
$r=\dim(L)$. 
 From theorem \nmb!{2.6} applied to the action of the Weyl group 
$W$ on $\Si$, it follows that the linear span of  
$\operatorname{grad}\si_1(x),\dots,\operatorname{grad}\si_n(x)$ 
equals $(T_x(W.x)^\bot)^{W_x}=\Si^{W_x}$, which is just the 
intersection of all walls containing $x$. 
For generic $x\in L$ we choose two different basis 
$\operatorname{grad}\si_{i_1}(x),\dots,\operatorname{grad}\si_{i_r}(x)$ 
and 
$\operatorname{grad}\si_{j_1}(x),\dots,\operatorname{grad}\si_{j_r}(x)$ 
of $L$, then the Gram matrix 
$$\align
B_{i_1,\dots,i_r}^{j_1,\dots,j_r}&= 
\pmatrix \langle d\si_{i_1}|d\si_{j_1}\rangle & \hdots & 
                    \langle d\si_{i_1}|d\si_{j_r}\rangle\\
                    \vdots &\ddots &\vdots\\
                 \langle d\si_{i_r}|d\si_{j_1}\rangle & \hdots &
                    \langle d\si_{i_r}|d\si_{j_r}\rangle
          \endpmatrix =\\
&=\pmatrix \frac{\partial \si_{i_1}}{\partial x_1}&\hdots& 
          \frac{\partial \si_{i_1}}{\partial x_n}\\ 
     \vdots&\ddots&\vdots\\ 
     \frac{\partial \si_{i_r}}{\partial x_1}&\hdots&
          \frac{\partial \si_{i_r}}{\partial x_n}\endpmatrix 
\pmatrix \frac{\partial \si_{j_1}}{\partial x_1}&\hdots& 
          \frac{\partial \si_{j_1}}{\partial x_n}\\ 
     \vdots&\ddots&\vdots\\ 
     \frac{\partial \si_{j_r}}{\partial x_1}&\hdots&
          \frac{\partial \si_{j_r}}{\partial x_n}\endpmatrix^\top. 
\endalign$$
has rank $r$ and thus for its determinant 
$$
\De_{i_1,\dots,i_r}^{j_1,\dots,j_r} = \langle  
     \operatorname{grad}\si_{i_1}(x)\wedge\dots
     \wedge\operatorname{grad}\si_{i_r}(x),
     \operatorname{grad}\si_{j_1}(x)\wedge\dots
     \wedge\operatorname{grad}\si_{j_r}(x) \rangle
$$
we have that 
$\De_{i_1,\dots,i_r}^{j_1,\dots,j_r}\o \bar c_1$ does not vanish 
identically since $\bar c_1$ hits generic points of $L$. 
Any minor of size $r$ which does not vanish identically along 
$\bar c_1$ is constructed by choosing two such bases, since for 
$\bar c_1(t)$ generic in $L$ the minor is not zero there exactly if 
the corresponding gradients form bases of $L$. 
But if we assume that $\bar c_1$ hits some smaller intersection $L_1$ 
of walls at times $t$ which are not isolated, then they cluster at a 
point $t_0$. 
By lemma \nmb!{4.3} there are real analytic orthogonal vector fields 
with the same wedge as the gradients, 
and one of them has to vanish along $L_1$. This 
one composed with $\bar c_1$ is then flat at $t_0$, thus the 
determinant is also. 
Hence each function
$\De_{i_1,\dots,i_r}^{j_1,\dots,j_r}\o \bar c_1$
has a zero of infinite order at 
$t=t_0$. This contradicts \nmb!{4.1}.(2), since the $r$ here is the 
same as there: Any larger minor vanishes on $L$ since any $s>r$ of 
the gradients are linearly dependent on $L$. 
So we may conclude the proof as in the real analytic case.
\qed\enddemo

\proclaim{\nmb.{4.3}. Lemma. Gram-Schmidt procedure}
Let $x_i:\Bbb R\to V$ be real analytic or smooth curves, defined near 
0, such that $t\mapsto x_1(t)\wedge \dots\wedge x_n(t)$ does not 
vanish of infinite order at 0.

Then there exist real analytic or smooth curves 
$\bar x_i:\Bbb R\to V$, defined near 0, such that 
$\langle \bar x_i,\bar x_j \rangle=0$ for $i\ne j$ and 
$x_1(t)\wedge \dots\wedge x_k(t) = 
\bar x_1(t)\wedge \dots\wedge \bar x_k(t)$, for all $k$ and all $t$.
\endproclaim

\demo{Proof}
We use induction on $n$. For $n=1$ we put $\bar x_1=x_1$.
So let $x_1,\dots, x_{n+1}$ be given, and by induction we may assume 
that $x_1,\dots, x_n$ are pairwise orthogonal. 
No $x_i$ is infinitely flat at 0 by assumption, so we may write 
$x_i(t)=t^{n_i}y_i(t)$ where $y_i(0)\ne 0$ and $y_1,\dots,y_n$ are 
still everywhere orthogonal.
Then
$$
\bar x_{n+1}(t) := x_{n+1}(t) - \sum_{i=1}^{n}
\frac{\langle y_i(t),x_{n+1}(t) \rangle}{\langle y_i(t),y_i(t) \rangle}
y_i(t)
$$
solves the problem.
\qed\enddemo


\enddocument